\input amstex
\documentstyle{amsppt}
\NoRunningHeads
\magnification=1200
\NoBlackBoxes
\topmatter

\title 
Nonlinear Stokes Phenomena analytic classification invariants via generic perturbation
\endtitle
\author A.A.Glutsuk
\endauthor
\abstract
For a generic deformation of a two-dimensional holomorphic 
vector field with an 
elementary degenerate singular point (saddle-node) we express   
the  Martinet - Ramis orbital analytic classification invariants of 
the nonperturbed field in terms of the limit transitions  
between the linearizing charts of singularities of the perturbed field. In the case, when the multiplicity of the singular point of the nonperturbed field is 
equal to 2, 
we show that the Martinet - Ramis invariants are limits of transition 
functions that compare appropriately normalized canonic first integrals of 
the perturbed vector field in the linearizing charts. We prove a generalization of this statement for higher multiplicity singularities. 
For a generic deformation of higher-dimensional 
holomorphic vector field with saddle-node singularity we show that appropriate sectorial central 
manifolds of the nonperturbed vector field are limits of appropriate separatrices of singularities of the perturbed field. 
We prove the analogues of the two first results  
for Ecalle-Voronin analytic classification invariants of one-dimensional conformal maps tangent to identity.
\endabstract

\thanks Research supported by grant of Pro-Mathematica foundation of 
French Mathematical Society, by CRDF grant RM1-229a,  by INTAS 
grant 93-0570-ext, by Russian Foundation for Basic Research (RFFI) grant 98-01-00455,  by Russian State Fellowship for young scientists of 
Russian Academy of Sciences.\endthanks
\address Permanent address:\ \ \ \  
Independent \ University\  of \ Moscow,\ \ \  Bolshoi\ \  Vlasievskii pereulok, 11,
121002, Moscow, Russia\endaddress  
\address Present address: Max-Planck Institut f\"ur Mathematik, Vivatsgasse, 7, 53111, Bonn, Germany\endaddress
 
\endtopmatter
\document 
\define\G{\Gamma}
\head 1. Introduction \endhead

 \define\g{\Gamma}
 \define\dl{\delta}

 \define\he{H_{\var}}
\redefine\ge{\Gamma_{\var}}
\define\qe{q_{\var}}
\redefine\ts{\widetilde{S}}
\redefine\ae{\a(\var)}
\define\cd{\Bbb C^2}
\redefine\c{\Bbb C}
\define\hh{\widehat H}
\redefine\tz{\tilde z}
\redefine\ts{\widetilde S}
\redefine\Im{\operatorname{Im}}
\definition{Definition} Say that a singular point of a holomorphic vector field 
in a domain in $\cd$ is {\it saddle-node}, if the correspondent linearization operator 
has exactly one zero eigenvalue.
\enddefinition

\definition{Definition} Two holomorphic vector fields are said to be 
{\it orbitally analytically equivalent}, 
if there exists a biholomorphic diffeomorphism that maps phase curves of 
one vector field to those of the other.
\enddefinition

\remark{Remark 1} Any germ of a holomorphic vector field in $(\cd,0)$ 
with saddle-node singularity at 0 is orbitally analytically equivalent to 
a vector field of the type

$$\cases
& \dot z=z+O(|z|^2+|t|^{k+1})\\
& \dot t=t^{k+1}\endcases \tag1$$
(see \cite{1}).
\endremark

One can ask the question: is it possible to transform a vector field 
(1) to a field that defines a differential equation with separable variables, 
more precisely, is is true that (1) is orbitally analytically equivalent to a direct sum of one-dimensional vector fields? Generically, the answer is "no" 
\cite{1}. At the same time, there always exists a formal Taylor series 
$\tz=\hh(z,t)$ such that the correspondent formal change of the coordinates 
$(z,t)\mapsto (\tz,t)$ is invertible and transforms (1) to a nonzero holomorphic 
function multiple of a unique vector field of the type 

$$\hskip-4cm(1)_n\hskip4cm\cases
 \dot{\tz}=\tz\\ 
 \dot t=t^{k+1}(1+\lambda t^k)^{-1},
\endcases$$  
$\lambda\in\c$, which is called the {\it formal normal form} of (1). 

Generically, the normalizing series $\hh$ diverges \cite{1}. At the same time, there 
exists a covering of a punctured neighborhood of zero in the $t$- plane by $2k$ 
sectors $S_j$, $j=0,\dots2k-1$, with the vertex at 0 having the following property: 
there exists a holomorphic change $\tz=H_j(z,t)$ of the variable $z$ in a 
neighborhood of 
zero in the domain $\ts_j=\c\times S_j$ that transforms phase curves of (1) to 
those of the formal normal form $(1)_n$ and has asymptotic Taylor series at 0 
that coincides with $\hh$ (\cite{1}, \cite{2}). Each sector $S_j$ contains the radial ray 
with the argument $\frac{\pi(1+2j)}{2k}$ from the set $\{ t^k\in i\Bbb R\}$
 and its closure does not contain an additional ray from this set.

The nontriviality of transitions between these normalizing charts yields  
the obstruction for orbital analytic equivalence between (1) and $(1)_n$. 
This phenomena is called the nonlinear Stokes phenomena. 
The transitions 
present a complete system of orbital analytic equivalence invariants of 
saddle-node singularities of vector fields with a given formal normal form, which 
are called Martinet - Ramis invariants \cite{1,2} and defined as follows. The formal 
normal form $(1)_n$ has the canonic first integral 
$$ I(\tz,t)=\tz t^{-\lambda}e^{\frac1{kt^k}}.$$
The first integral $I$ together with the sectorial normalizing coordinate changes $H_j$ 
define canonic first integrals $I_j=I\circ H_j$ of (1) in the domains $\ts_j$. (The sectors are numerated in the 
counterclockwise direction. We identify $S_{2k}$ with $S_0$, and  
put also $I_{2k}=I\circ H_0$. In the definitions of all the integrals $I_j$ 
we consider that for any pair of neighbor sectorial domains $\ts_j$ and 
$\ts_{j+1}$, $j<2k$, the 
branch of the (multivalued) integral $I$ in the latter  is obtained 
from that in the former by counterclockwise analytic continuation.)  Let   
$\phi_j$ be the transition functions that compare the first integrals $I_j$  
and $I_{j+1}$ over the intersection component of the correspondent neighbor 
sectors $S_j$ and $S_{j+1}$: $I_{j+1}=\phi_j\circ I_j$. 

\remark{Remark 2} The functions $\phi_j(\tau)$ are holomorphic in $\c$ and have 
the type $\phi_j(\tau)=\tau+a_j$, if 
$j$ is odd; for even $j$ they are holomorphic in a neighborhood of zero and have  unit derivative at 0: $\phi_j(\tau)=\tau+o(\tau)$, as $\tau\to0$. Two germs of vector 
fields in $(\cd,0)$ of type (1) are orbitally analytically equivalent, 
iff they have the same formal normal form and the correspondent sets of $2k$ transition functions 
$\phi_j$, $\tilde{\phi_j}$ ($j=0,\dots,2k-1$) can be obtained from each other by simultaneous 
conjugation by multiplication by constant, $\phi_j(\tau)= 
c\tilde{\phi_j}(c^{-1}\tau)$, and cyclic permutation. 

The equivalence 
classes of the transition function sets with respect to the last operations are 
called Martinet - Ramis invariants.  
\endremark

\example{Example 1} Let $k=1$. Then the above covering of a punctured 
neighborhood of zero consists of two sectors: $S_0$ and $S_1$ (Fig.1). 
The sector $S_0$ contains the positive imaginary semiaxis, and its  closure 
does not contain the negative imaginary semiaxis; the other sector $S_1$ 
possesses the same properties with the interchange of "negative" and "positive".  
In this case we have two transition functions $\phi_0$ 
and $\phi_1$ correspondent to the two components of the intersection of 
the sectors.
\endexample
\hskip-2cm {\bf Fig.1}

In the present paper we consider a generic continuous deformation

\redefine\le{(1)_{\var}}
\define\a{\alpha}
\define\var{\varepsilon}
\define\pr{\prod_{i=0}^{k}(t-\a_i(\var))}
\define\ije{I_{i,\var}}
\define\aje{\a_j(\var)}

$$\hskip-2cm(1)_{\var}\hskip2cm\cases 
\dot z=z(1+q(z,t,\var))+g(t,\var)\pr\\
\dot t=\pr\endcases$$
of a field (1) in the class of holomorphic vector fields, $q(0,0,0)=0$. In particular, genericity 
means that the degenerate singularity 0 of the nonperturbed vector field 
(1) is split into $k+1$ nondegenerate linearizable singularities 
$(0,\a_i(\var))$ of the perturbed field, so that in a neighborhood of each of the 
latters there exists a local biholomorphic change of variables that transforms
$\le$ to a linear field 

$$\cases\
& \dot z=z\\ 
& \dot t=\mu_i(\var)t,
\endcases\tag2$$
where $\mu_i(\var)\to 0$, as $\var\to0$. The field (2) has the 
canonic first integral $zt^{-\mu_i(\var)}$. The latter together with 
the linearizing change of variables defines a canonic first integral $\ije$
of $\le$ in a neighborhood of the singular point $(0,\a_i(\var))$. The integral 
family $I_{i,\var}$ is uniquely defined up to multiplication by constant family.  

The main result of the paper is the following. In the case, when $k=1$, for appropriately normalized  
canonic first integrals $\ije$ we show that the transition 
functions that compare them in the intersection components of the correspondent 
linearization domains tend to the transition functions $\phi_i$ that define 
Martinet - Ramis invariants of the nonperturbed field (1), as $\var\to0$. We 
prove a generalization of this statement for larger $k$. 
In fact, for any $k\in\Bbb N$ we prove, that appropriate single-valued branches of canonic first integrals of the perturbed 
field tend to appropriate sectorial canonic first integrals of the nonperturbed field. This result is stated in Subsection 2.3 and proved in Subsections 3.1-3.3 and 3.5.

In Section 4 we prove its higher-dimensional "central manifold" analogue. 

 A similar statement on the expression of Ecalle - Voronin moduli for 
 analytic classification of germs of conformal maps $z\mapsto z+z^{k+1}+\dots$ 
via perturbation is stated in Subsection 2.4 and proved in Subsection 3.4. After obtaining a proof for some class of deformations, 
the author became aware that for deformations analytic in the parameter 
this statement was known to specialists (as A.Douady)  
and in the simplest case $k=1$ a proof was obtained by J.Martinet \cite{4}. The proof in the general case was not published. 
After being acquainted with the known proofs, the author  
found that the method he used is quite different from that of the latters.  
It is similar to that of the previous statement on the Martinet-Ramis invariants. 
 
Earlier a similar problem to obtain analytic classification invariants 
of linear differential equation with nonresonant irregular singularity in terms 
of its deformation 
was studied by R.Garnier \cite{3}, C.Zhang \cite{5}, J.-P.Ramis \cite{12}, 
A.Duval \cite{13}, the author \cite{6} and others (see the references in \cite{6}). 
In 1919 for a special class of such equations R. Garnier had obtained analytic 
classification invariants by studying some their particular deformations \cite{3}.  In 1984 V.I.Arnold conjectured that Stokes operators of irregular singularity of linear differential equation can be expressed in terms of limit 
monodromy of its deformation. This conjecture was also independently stated by J.-P.Ramis in 1988 who proved that for classical confluenting family of hypergeometric equations Stokes operators of the nonperturbed equation are limits of transition operators between appropriately normalized monodromy eigenbases of the perturbed equation \cite{12}. In 1991 A.Duval \cite{13} 
proved the same statement for biconfluenting family of hypergeometric equations.
In 1994 C.Zhang had obtained the expression of Garnier invariants via 
Stokes operators \cite{5}. In fact, Garnier and Zhang implicitly proved the previous statement for the particular deformations of the equations considered by Garnier. In 1997 this statement was proved by the author for generic deformation of any irregular nonresonant singularity \cite{6}.

\head 2. Main results \endhead

To state the main results, we recall some definitions, the Sectorial 
Normalization Theorem for saddle-node singularity of holomorphic vector field 
mentioned in Section 1, and the linearization Theorem for generic vector field with nondegenerate singularity. 

\subhead 2.1. Saddle-node vector fields. Sectorial Normalization Theorem 
\endsubhead

\definition{Definition} Let $k\in \Bbb N$. An {\it imaginary dividing ray} is a 
radial ray from the set $\{ t^k\in i\Bbb R\}$ in the complex plane with the coordinate 
$t$. A radial sector $S\subset\Bbb C$ is said to be 
{\it good}, if it contains a unique imaginary dividing ray 
and its closure does not contain an additional such ray.   
It is said to be $j$- good, if the correspondent ray has the argument  
$\frac{\pi(1+2j)}{2k}$ (the latter will be referred to, as $r_j$).
\enddefinition

\example{Example 2} Let $k=1$. Then the imaginary dividing rays are the 
imaginary semiaxes. The sectors from Example 1 are good.
\endexample

\define\const{\operatorname{const}}

\remark{Remark 3} The imaginary dividing rays are exactly the rays, 
where the asymptotics of the first integral $I$ of a vector field $(1)_n$ 
(Section 1), as  $ t\to0$, $z=\const$, is changed from zero to 
infinity. 
\endremark
 
\definition{Definition} Let $(z,t)$ be a coordinate system in $\cd$, $U$ be a neighborhood of zero. For 
a subset $S$ in the $t$- plane define $\ts=(\c\times S)\cap U$. 
\enddefinition

\proclaim{Theorem 1}(Sectorial Normalization Theorem \cite{1,2}.) For any 
vector field (1) there exists a neighborhood $U$ of zero in the phase space 
such that for any $j$- good sector $S_j$ there exists a change 
$\tz=H_j(z,t)$ of the 
variable $z$ holomorphic in the domain $\ts_j$ and 
$C^{\infty}$- smooth in its closure that transforms phase 
curves of (1) to those of its formal normal form $(1)_n$  (see Section 1) and has asymptotic Taylor series at 0 coinciding with the normalizing series $\hh$ 
from the beginning of Section 1. 
\endproclaim

\definition{Definition} In the conditions of Theorem 1 the first integral 
of the field (1) in $\ts_j$ obtained from the canonic first 
integral $I$ of the formal normal form by sectorial 
normalizing chart is called the sectorial {\it canonic  first integral}. 
\enddefinition

\remark{Remark 4} The sectorial canonic first integral is uniquely 
defined up to multiplication by constant \cite{1}.
\endremark

\subhead 2.2. Linearization Theorem \endsubhead

\definition{Definition}\ \  \ \ \ The {\it characteristic number} \ \ of\ a holomorphic vector 
field \ \ in two-dimensional complex domain at its singular point is the ratio of the eigenvalues of the 
correspondent linearization operator.
\enddefinition

\definition{Definition} A singular point of a holomorphic vector field 
in a two-dimensional complex domain is said to be {\it linearizable}, if the field is locally orbitally analytically equivalent to its linear part at the 
singularity in its neighborhood.
 \enddefinition

\proclaim{Theorem 2 \cite{7}} Let a holomorphic vector field in $(\cd,0)$ 
have singular point at 0 with finite nonreal characteristic number. Then the singular point is linearizable.
\endproclaim

\definition{Definition} The {\it canonic first integral} of a linear vector field 
$$\cases\
& \dot z=\lambda z\\ 
& \dot t=\mu t,
\endcases$$
with $|\lambda|>|\mu|$  is  $zt^{-\frac{\lambda}{\mu}}$. 
\enddefinition

\definition{Definition} Let a holomorphic vector field have a linearizable   singular 
point with linear part as in the previous Definition. The {\it canonic first  integral} of the field at the singular point is the first integral obtained from  the canonic first integral of the linear part by linearizing chart. 
\enddefinition

\remark{Remark 5} In the conditions of the previous Definition the canonic 
first integral is uniquely defined up to multiplication by constant.
\endremark

\subhead 2.3. Main results. Martinet - Ramis invariants are limits of 
transition functions between canonic integrals of generic perturbation \endsubhead

We consider a continuous deformation $\le$ of a vector field (1) from Section 1 
that in particular splits the degenerate singular point 0 into $k+1$ 
nondegenerate singularities $(0,\a_i(\var))$ of the perturbed field, i.e., $\a_i(\var)\neq 
\a_l(\var)$ for $i\neq l$, $\var\neq0$. For a generic deformation $\le$ we prove the 
statement from the title of the Subsection.

\remark{Remark 6} We restrict ourselves by considering only deformations of the type 
$\le$ without loss of generality. Namely, in Subsection 3.7 of Section 3 we 
prove the following

\proclaim{Lemma 1} Any continuous deformation of any vector field (1) 
in the class of holomorphic vector fields (without the condition that 
all the singularities of the perturbed field are nondegenerate, the only 
requirement is continuity) is orbitally analytically equivalent 
to a family of the type $\le$ via continuous family of changes of variables 
defined for all parameter values small enough. 
\endproclaim
\endremark

\subhead 2.3.A. $k=1$\endsubhead    
Firstly let us state the main result in the case, when $k=1$ (i.e., the perturbed field has two singularities). To do this, 
let us introduce the following

\definition{Definition} A family $\le$ as at the beginning of the Subsection 
is said to be {\it nondegenerate}, if the line in the $t$- plane passing through the singularity pair of 
the perturbed field intersects the real axis by angle bounded away from zero 
uniformly in $\var$ small enough. 
\enddefinition

\remark{Remark 7} A family $\le$ is nondegenerate, iff the characteristic numbers of the singular points of the perturbed vector field have arguments bounded 
away from $\pi\Bbb Z$ uniformly in $\var$ small enough. In particular, in this 
case the singularities of the perturbed field are linearizable (Theorem 2). For small parameter 
values the singularities of the perturbed field satisfy the conditions of 
the Definition of the canonic first integral at a linearizable singular point 
from the previous Subsection. Thus, the canonic first integrals 
of the perturbed vector field at its singularities are well-defined for any $\var$ small enough. 
\endremark

We consider nondegenerate families $\le$.  
Without loss of generality we consider that $\a_0=-\a_1$.  

\remark{Remark 8} Let a family $\le$ be as in the last item. Then for 
small parameter values $\Im \a_j(\var)$ have constant signs. 
The arguments $\arg\a_j$ are bounded away from $\pi\Bbb Z$. We consider that 
$\Im \a_0(\var)>0$, $\Im\a_1(\var)<0$. 
\endremark

To the family $\a_j$ we put into correspondence the sector $S_j\supset\{(-1)^j
\operatorname{Im}t>0\}$ from 
Example 1 (Fig.1, 2a), which contains $\a_j(\var)$ for all $\var$. We show that the branch of appropriately normalized 
(multivalued) canonic first 
integral of the perturbed field at its singular point $(0,\a_j(\var))$ 
converges to the sectorial canonic first integral of the nonperturbed vector 
field in $\ts_j$. 
\redefine\D{\Delta}
\redefine\d{\delta}

\define\ie{I_{\var}}

\define\ov{\Omega_{\var}}
\define\Lim{\operatorname{Lim}}

\definition{Definition} Let $r>0$, $S$ be a radial sector in complex line 
with the coordinate $t$. Define $S^r=S\cap\{|t|<r\}$.
\enddefinition

\proclaim{Theorem 3} Let $k=1$, $\le$ be a nondegenerate vector field family, $\a(\var)$, $-\a(\var)$ be the 
continuous families of the $t$- coordinates of its singular points, $S$ 
be the sector correspondent to $\a(\var)$  from the item following Remark 8. There exist an $r>0$, neighborhoods 
$U_z=\{|z|<2\d\}$, $U_t=\{|t|<\d\}$ of zero in the $z$- and $t$- axes respectively and a family $\ov$ of simply 
 connected subdomains of the disc $U_t$ that contain $\a(\var)$,  
do not contain $-\a(\var)$ and possess the following properties: 
 
 1)\footnote{A family $V_{\var}$ of planar domains is said to be convergent to 
 a domain $V$,  if both the maximal distance between a point 
 of the boundary $\partial V_{\var}$  and the whole boundary $\partial V$, 
 and that with the interchange of $V$ and $V_{\var}$ tend to 0}
the connected component of the intersection $\ov\cap (S^r\setminus[0,-\a(\var)])$ containing $\a(\var)$ tends to $S^r$, as $\var\to0$.

 2) Put $\ov'=\ov\setminus[-\a(\var),\a(\var)]$. The canonic first integral $\ie$ of the perturbed vector field $\le$ 
 at its singularity $(0,\a(\var))$ is a multivalued holomorphic function 
 in the domain $\widetilde{\ov}=U_z\times\ov$ with branching at the line $t=\a(\var)$. 
 It is single-valued in its subdomain $\widetilde{\ov'}$. The restriction to the latter of appropriately normalized integral $I_\var$  converges\footnote{In the conditions of the previous footnote a family of functions holomorphic in $V_\var$ depending continuously on the same parameter $\var$ is said to be convergent (in $V$), if it converges uniformly in compact subsets of $V$}
 to the canonic 
 sectorial integral of the nonperturbed vector field in the domain $\widetilde{S^r}=
U_z\times S^r$. 
 \endproclaim

Theorem 3 is proved in Subsections 3.1-3.3 and 3.5.

\remark{Remark 9} The domain $\ov$ correspondent to $\a=\a_0$ is depicted at 
Fig.4a. It will be shown that it converges to the domain $\Omega$ depicted at 
Fig.4b bounded by a cardioid-like curve having inward cusp at 0 with two 
distinct tangent rays that bound a sector disjoint from both $S$ and $\Omega$. 
\endremark
 
 \proclaim{Corollary 1} In the conditions of Theorem 3 let $\a_j$ be singularity $t$- coordinate families of the vector fields $\le$, $S_j^r$, $\ov'=\ov'(j)$, be the correspondent sectors and domains from Theorem 3, $C_0$ be a connected component of the intersection $S_0^r\cap S_1^r$. Let $\phi$ be the correspondent Martinet-Ramis transition function between the sectorial canonic 
 integrals of the nonperturbed field from Theorem 1 over $C_0$. There exists a family $C_\var$ of connected components of the intersections $S_0^r\cap S_1^r\cap\ov'(0)\cap\ov'(1)$ that converges to $C_0$ and possesses the following property.  The transition functions between appropriately normalized canonic first integrals of the 
 perturbed field in $\widetilde{C_\var}$ are well-defined and holomorphic in a  domain (depending on $\var$) that converges to the definition domain of the function $\phi$ (either neighborhood of zero, or $\Bbb C$, see Remark 2). These transition functions converge to $\phi$ (footnote 2). \endproclaim 

\hskip-2cm {\bf Fig.2a,b}

\subhead 2.3.B. $k\geq2$ \endsubhead 
Now let us state the main result in the case,  when $k\geq2$. 
To do this, let us extend the Definition of nondegenerate deformation. 

We consider  families $\le$ with the following properties:
 the polynomial family $p(t,\var)=\pr$ is differentiable in $\var$ at 
$\var=0$, $p'_{\var}(0,0)\neq0$.
Without loss of generality we consider that $\sum\a_j=0$.

\remark{Remark 10} Let $p(t,\var)$ be a family as in the last item. 
Then its root polygon is asymptotically regular. This means that 
its radial homothethy image with the diameter 1 tends to a regular polygon (with the center 
at 0), as $\var\to0$. The latter will be referred to, as $\D$. Its vertex that 
 is the limit of the homothethy images of the points $\a_j(\var)$ will be 
 denoted by $A_j$. (We consider that $\a_j$ (and hence, $A_j$) are numerated in the 
 counterclockwise direction.)
\endremark

\definition{Definition} Let $k\in\Bbb N$. A {\it real (imaginary) dividing ray} or line  
in the complex plane with the coordinate $t$ is a radial ray or line where 
$t^k\in\Bbb R$ (respectively, $t^k\in i\Bbb R$).
\enddefinition

\definition{Definition} A family $\le$ as at the beginning of Subsection B 
is said to be 
{\it nondegenerate}, if no bissectrix of the correspondent limit regular polygon $\D$ 
from the previous Remark lies in a real dividing line.
\enddefinition

\remark{Remark 11} For a family $\le$ as at the beginning of the 
Subsection, the nondegeneracy is equivalent to the characteristic number  
property of the singularities of the perturbed vector field from Remark 7. It 
is also equivalent to the condition that for any 
$i=0,\dots,k$ there is a (unique) imaginary dividing ray 
$r_{j_i}$ that has angle less than $\frac{\pi}{2k}$ with the radial ray 
of $A_i$ (so, $r_{j_i}$ is the imaginary dividing ray closest to the radial ray). 
For a nondegenerate family $\le$, to the singularity coordinate family $\a_i$ we put into correspondence a $j_i$- good
sector $S_{j_i}$ that contains $r_{j_i}$ and the radial ray such that the latter has angles greater than $\frac{\pi}{2k}$ with the boundary rays of $S_{j_i}$  (see Fig.2b for $k=2$). Without loss of generality, everywhere below we consider that $j_0$ is equal to either 0, or 1 and the sequence $\{j_i\}_{i=0,\dots,k}$ is strictly increasing). One can achieve this by applying appropriate linear change of the variable $t$ of the type $t\mapsto e^{ 2i\frac{\pi}k}t$. We put  $\a_{k+1}=\a_0$, $j_{k+1}=j_0+2k$, $S_{j_{k+1}}=S_{j_0}$. 
\endremark
 
\remark{Remark 12} In the conditions of the previous Remark the sectors 
correspondent to the singularity families do not cover punctured neighborhood of zero: the total number of the imaginary dividing rays is $2k$, and that of the sectors is $k+1$ (each of the latters contains exactly one imaginary dividing ray). On the other hand, for each pair of sectors $S_{j_i}$ and 
$S_{j_{i+1}}$, $i\leq k$, correspondent to neighbor singularity coordinate families $j_{i+1}-j_i\leq2$, i.e., either the sectors intersect each other ($j_{i+1}-j_i=1$), or they are intersected by the unique intermediate sector $S_{j_i+1}$ of the covering ($j_{i+1}-j_i=2$). Indeed, the angle between each neighbor pair of the $A_i$ radial rays from the previous Remark is equal to 
$\frac{2\pi}{k+1}$. Therefore, the angle between the correspondent imaginary dividing rays is less than $\frac{2\pi}{k+1}+\frac{\pi}k<3\frac{\pi}k$ (by definition), and hence, not  greater than $2\frac{\pi}k$ (since this angle is a multiple of $\frac{\pi}k$ by definition).  There are exactly two intersected sector pairs 
$(S_{j_i}, \ S_{j_{i+1}})$ (i.e., with $j_{i+1}=j_i+1$). 
\endremark

\proclaim{Theorem 4} Let $\le$ be a nondegenerate vector field family, 
$\a=\a_i$, $S=S_{j_i}$ be respectively its singularity coordinate family and 
the correspondent sector from Remark 11. There exist a neighborhood $U=U_z\times U_t$ of zero in the phase space and a family of domains 
$\ov=\ov(i)\subset U_t$ such that the triple $(\a,\ov, S)$ satisfies the statements of Theorem 3 with the change of $[0,-\a(\var)]$ in its statement 1) 
to $\bigcup_{\a_s\neq\a}[0,\a_s(\var)]$ and $[-\a(\var),\a(\var)]$ in its statement 2)  to $\bigcup_s[0,\a_s(\var)]$. 
 \endproclaim
 
 Theorem 4 is proved in Subsections 3.1-3.3 and 3.5.  

Theorem 4 admits a Corollary analogous to Corollary 1 for the transition functions $\phi_{j_l,j_{l+1}}$ between the canonic first integrals of the nonperturbed field correspondent to the disjoint sector pairs $S_{j_l}$, $S_{j_{l+1}}$ (i.e., with $j_{l+1}=j_l+2$, the number of these pairs is equal to $k-1$, see the previous Remark). These transition functions are defined in the next item. All the Martinet-Ramis transition functions 
(different from those correspondent to the (two) intersected sector pairs 
$(S_{j_l},\ S_{j_{l+1}})$) are  expressed in terms of the functions $\phi_{j_l,j_{l+1}}$ at the end of the next item. 
 
\remark{Remark 13} Let (1) and $\bigcup_{j=0}^{2k-1}S_j$ be respectively a vector field and a covering from Section 1. Let $S_l$, $S_{l+1}$ and $S_{l+2}$ be a triple of consequent sectors from the covering, $l\leq 2k-1$ (we consider that $S_{2k}=S_0$, $S_{2k+1}=S_1$). Let $I_l$, $I_{l+1}$, $I_{l+2}$ be the correspondent canonic integrals (defined consequently as in the item preceding Remark 2 in Section 1). Let $\phi_l$ and $\phi_{l+1}$ be the correspondent transition functions from Section 1. (Then $\phi_{2k}(\tau)=e^{-2\pi i\lambda}\phi_0(e^{2\pi i\lambda}\tau)$, by definition and since the analytic 
continuation by counterclockwise going around zero transforms the canonic integral $I$ of the formal normal form to $e^{-2\pi i\lambda}I$.) In the case, when $l$ is even, the zero level curves of the integrals $I_l$ and $I_{l+1}$ coincide over the intersection component of the correspondent sectors, i.e., continue analytically each other (Remark 2). In particular, $I_l$ extends analytically counterclockwise to a neighborhood of this conjoint zero level curve  
in $\tilde S_{l+2}\cap \tilde S_{l+1}$, and the transition function $\phi_{l, l+2}(\tau)$ 
between the restrictions of $I_l$ and $I_{l+2}$ to this neighborhood   ($I_{l+2}=\phi_{l,l+2}\circ I_l$) is well-defined  
and univalent in a neighborhood of zero. This transition function is equal to the composition 
$\phi_{l+1}\circ\phi_l$. {\it The transition functions 
$\phi_l$ and $\phi_{l+1}$ are expressed in terms of the composition function $\phi_{l,l+2}$}  as follows: 
$$\phi_{l+1}(\tau)=\tau+\phi_{l,l+2}(0),\ \ \ \ \phi_l(\tau)=\phi_{l,l+2}(\tau)-
\phi_{l,l+2}(0).$$ 
Analogously, in the case, when $l$ is odd, the inverse transition function 
$\phi_{l,l+2}^{-1}$ is well-defined and univalent in a neighborhood of zero and 
$$\phi_{l,l+2}^{-1}=\phi_l^{-1}\circ\phi_{l+1}^{-1},\ \ \  \phi_l(\tau)=\tau-\phi_{l,l+2}^{-1}(0),\ \ \  \phi_{l+1}(\tau)=\phi_{l,l+2}(\tau+\phi_{l,l+2}^{-1}(0)).$$
\endremark

\proclaim{Corollary 2} In the conditions of Theorem 4 let $\a_l$, $\a_{l+1}$ be neighbor singularity coordinate families of $\le$, $S_{j_l}^r$,  $S_{j_{l+1}}^r$, $\ov(l)$, $\ov(l+1)$ be the correspondent sectors and domain families. Let $\ov'(s)=\ov(s)\setminus\bigcup_{l=0}^k[0,\a_l(\var)]$, $s=l,l+1$.     In the case, when the  sectors are intersected ($j_{l+1}-j_l=1$), 
the tuple $S_{j_l}^r$, $S_{j_{l+1}}^r$, $\ov'(l)$, $\ov'(l+1)$ satisfies the statements of Corollary 1. Otherwise, in the case, when the sectors are disjoint   ($j_{l+1}-j_l=2$, 
see Remark 12) there exists a family $C_\var$ of subdomains in the union  $S_{j_l}^r\cup S_{j_l+1}^r$, $\a(\var)\in C_\var$, converging to the latter with the following properties:

1) Case, when  $j_l$ is even. The zero level curve of the canonic first integral $I_{l,\var}$ of the perturbed field at its singular point $(0,\a_l(\var))$ 
extends up to an analytic curve containing   
the graph $z=\qe(t)$ of a function $\qe(t)$ holomorphic in $C_\var$.   
The integral $I_{l,\var}$ extends analytically to a neighborhood of this graph. (This neighborhood intersects $\widetilde{\ov'}(l+1)$, which is the definition domain  of the canonic  integral $I_{l+1,\var}$, whenever $\var$ is small 
enough.) Let $\phi_{l,\var}$, $I_{l+1,\var}=\phi_{l,\var}\circ 
I_{l,\var}$, be the transition function between the extended integral $I_{l,\var}$ and the integral $I_{l+1,\var}$ in 
the intersection of the correspondent definition domains.  For appropriately normalized canonic integrals the function $\phi_{l,\var}$ is well-defined and holomorphic in a neighborhood of zero independent on $\var$ (for all $\var$ small enough) and tends to the composition $\phi_{j_l+1}\circ\phi_{j_l}$ of the Martinet-Ramis transition functions of the nonperturbed field (defined in Remark 2), as $\var\to0$.

2) Case, when $j_l$ is odd. The same statements are valid with respect to the 
clockwise analytic extension of the integral $I_{j_{l+1},\var}$ and its zero level curve (from $\widetilde{\ov'}(l+1)$ to a subdomain of $\widetilde{\ov'}(l)$), the 
inverse transition function $\phi_{l,\var}^{-1}$ and the composition $\phi_{j_l}^{-1}\circ\phi_{j_{l+1}}^{-1}$ (the function $\phi_{2k}$ is defined in the previous Remark). 
\endproclaim

 \subhead 2.4. Analogues of main results for germs of conformal maps with identity linear 
 part\endsubhead
 
 \define\la{\lambda}
 \define\vl_{v_{\lambda}}
 
 In this Subsection we state an analogue of 
 Theorems 3 and 4 for germs of one-dimensional conformal maps in 
 $(\c,0)$ with the fixed point 0 of unit multiplier. 
 To do this, let us introduce some definitions 
 and recall the analogue of nonlinear Stokes phenomena for such germs.
  
  \remark{Remark 14} Any germ as in the last item can be transformed to a one 
 of the type $t\to t+2\pi it^{k+1}(1+O(t))$, as $t\to0$, by linear change of the 
 coordinate.  
\endremark

\subhead 2.4.A.\ \  Nonlinear \ \ Stokes phenomena\ \  for\ \  one-dimensional conformal maps
\endsubhead

 \definition{Definition} Two germs of conformal mappings $(\c,0)\to(\c,0)$ are 
 said to be {\it analytically (formally) equivalent}, if there exists a germ of conformal 
 diffeomorphism (an invertible formal Taylor series) that conjugates them. 
 \enddefinition

\definition{Definition} Let $f$ be a conformal map defined in a domain $U$, 
$D\subset U$ be its subdomain, $v$ 
be a holomorphic vector field in $U$ such that all its time $s$ flow maps 
$g^s_v$ with $0<s\leq1$ map $D$ to $U$. The vector field $v$ is 
said to be the {\it generator} of $f$, if its unit time flow map is $f$:\ 
$g^1_v|_{D}=f$.
\enddefinition

\redefine\vl{v_{\lambda}}

The analytic classification of germs as in the previous Remark was obtained 
independently by J.Ecalle \cite{8} and S.M.Voronin \cite{9}. Namely, 
any such germ $f$ has a formal generator, more precisely, 
is formally equivalent to  the unit time flow map germ of a unique 
vector field of the type $\vl(t)=2\pi it^{k+1}(1+\la t^k)^{-1}$. This map  
is called the formal normal form of $f$. Generically the 
normalizing Taylor series diverges. At the same time there exists a covering of 
a punctured neighborhood of zero by $2k$ sectors $S_j$, $j=0,\dots,2k-1$, such 
that in each sector there exists a conformal diffeomorphism that conjugates $f$ 
to its formal normal form with the following properties:

\proclaim{Theorem 5 \cite{1,8,9}}  Let $f$ be a germ as in Remark 14. Then for any $j$- good sector $S_j$ (see the Definition in  
Subsection 2.1) there exists an $r>0$ and a conformal map $h_j:S_j^r\to \c$ that conjugates $f$ to its formal normal form with the following properties: 

1) the map $h_j$ is $C^{\infty}$ in the closure of the sector $S_j^r$, $h_j(0)=0$,  $h_j'(0)=1$; its asymptotic 
Taylor series at 0 conjugates $f$ to its formal normal form. 

2) the maps $h_j$ correspondent to different $j$ may be chosen to have 
common asymptotic Taylor series at 0.  
\endproclaim

\definition{Definition} Let $f$, $S_j^r$, $h_j$ be as in Theorem 5. The generator $v_j$ of the restriction of the map $f$ to $S_j^r$ obtained from the correspondent vector field $\vl$ by the normalizing chart $h_j$ is called the {\it canonic 
sectorial generator}.
\enddefinition

\remark{Remark 15} The sectorial canonic generator of a germ $f$ as in the 
previous Definition is the unique generator of $f$ that is holomorphic in $S_j^r$ 
and $o(t)$, as $t\to0$. Let us prove the uniqueness of 
such a generator. By Theorem 5, it suffices to do this for the formal normal 
form, which will be now denoted by $f$. Let $v$ be such a generator. Let us show that $v=v_{\lambda}$. For $r>0$ small enough there exists a domain $U_j$ containing the sector $S_j^r$ invariant with respect to the field $-\vl$ (in particular, 
$f^{-1}$- invariant) such that each backward $f$- orbit in $U_j$ converges to 0 in the asymptotic direction of the imaginary dividing ray in $S_j$ \cite{1} (in particular, it fits $S_j^r$ eventually). The space of these orbits is isomorphic to the double punctured Riemann sphere at 0 and $\infty$. The isomorphism is induced by the map $U_j\to\overline{\c}\setminus
\{0,\infty\}$ defined by the formula $t\mapsto \tau=t^{\la}e^{-\frac1{kt^k}}$ (its   
right-hand side is equal to $e^{2\pi iT}$, where $T$ is a complex time function 
of $\vl$).  Any generator of $f$ in $S_j^r$ that is $o(t)$, as $t\to0$ (in particular, $v$), induces a holomorphic vector field in the whole orbit 
space (let us denote this field correspondent to $v$ ($v_{\lambda}$) by $v'$ 
(respectively, $v'_{\lambda}$)). A straightforward calculation shows 
that $v'(\tau)\to0$,  
as $\tau$ tends to either 0, or $\infty$, so $v'$ continues  up to a holomorphic  vector field in the whole Riemann sphere that vanishes at the punctured points, 
and so, it is linear. It has the form $\dot{\tau}=2\pi im\tau$, $m\in\Bbb Z$: 
the correspondent unit time flow map should be identity, since that of the 
generator $v$ in $S_j^r$ is $f$.   
 The vector field $\vl'$ is 
 $\dot{\tau}=2\pi i\tau$, which follows from its definition. The two last statements  imply that $v=mv_{\lambda}$, so, $v$ is a generator of the $m$-th iteration of $f$. Therefore (since $f$ is not periodic), $m=1$ and $v=\vl$. 
 \endremark
 
\remark{Remark 16} In the conditions of the previous Definition let $\tau_j$ be 
the complex time function correspondent to the sectorial canonic generator in 
$S_j$ normalized so that $\tau_j-\tau_{j+1}\to0$, as $t\to0$ in the connected 
component of the intersection of the neighbor sectors $S_j$ and $S_{j+1}$, whenever 
$j\leq 2k-2$. Then $\tau_0-\tau_{2k-1}\to-\la$, as $t\to0$. (The possibility of 
such choice of time functions and the last statement are proved in 
\cite{1}.) We put $\tau_{2k}=\tau_0+\lambda$. Consider the transition 
functions $\psi_j(\tau)=\tau_{j+1}\circ \tau_j^{-1}(\tau)$ between the 
time charts in the intersection components of the correspondent 
sectors, $j=0,\dots,2k-1$. The functions $\psi_j(\tau)$  extend up to functions  holomorphic in appropriate half-planes $(-1)^j\Im\tau<c$, 
$c\in\Bbb R$, having the type 
$$\psi_j(\tau)=\tau+\sum_{(-1)^jl<0}c_je^{2\pi il\tau}.$$  
Two germs $f$ and $\tilde f$ as in Theorem 5 are analytically equivalent, iff  the correspondent 
transition function sets $\psi_j$ and $\widetilde{\psi_j}$ can be obtained one from 
the other by simultaneous conjugation by addition of constant,  
$\widetilde{\psi_j}(\tau)=\psi_j(\tau-c)+c$, and subsequent cyclic 
permutation of the functions \cite{1,8,9}. 
Thus, the equivalence classes of sets of transition functions with respect to 
the last operations are analytic classification invariants of the correspondent 
germs. They are called Ecalle - Voronin moduli \cite{1,8,9}.
\endremark

\subhead 2.4.B. Linearization of attracting and repelling fixed points. Canonical 
generators \endsubhead

\proclaim{Theorem 6 \cite{7}} Let a one-dimensional conformal map have a fixed point 
with multiplier of nonunit module. Then it is conjugated to its linear part 
by conformal diffeomorphism in some neighborhood of the fixed point.
The conjugating diffeomorphism is unique up to left composition with linear map. \endproclaim  

\definition{Definition} Let $\mu\in\c\setminus i\Bbb R_-$ be a number with a nonunit module. The canonic 
generator of the linear map $t\mapsto\mu t$ 
of the complex plane at its fixed point 0 is $\dot t=(\ln\mu) t$, 
where $\ln\mu$ 
is chosen so that $-\frac{\pi}2<\Im\ln\mu<\frac{3\pi}2$. The canonic generator of a nonlinear map at its fixed point with a multiplier as above is its generator 
obtained from the canonic generator of its linearization by the conjugating 
diffeomorphism from Theorem 6.
\enddefinition

\define\fe{f_{\var}}

\subhead 2.4.C. Ecalle - Voronin moduli via generators of perturbation\endsubhead

We consider a continuous deformation 
$$ \fe(t)=t+2\pi i(1+q(t,\var))\pr$$ 
of a conformal map $f$ as in Remark 14, $q(0,0)=0$. We show that for a generic 
deformation $\fe$ that in particular splits the degenerate fixed point of the 
nonperturbed map to $k+1$ nondegenerate linearizable fixed points $\a_i(\var)$ of the perturbed map appropriate 
sectorial canonic generators of the nonperturbed map are limits of 
canonic generators of the perturbed map.

\remark{Remark 17} Any continuous deformation of a map $f$ as in Remark 14 
in the class of conformal maps has the type $\fe$.
\endremark

\definition{Definition} A family $\fe$ as at the beginning of the Subsection is said to be {\it nondegenerate}, if 
the product $p(t,\var)=\pr$ in its formula satisfies the condition of one of the  two  Definitions 
of nondegenerate family from Subsection 2.3 (one distinguishes the cases  
$k=1$ and $k\geq2$). 
\enddefinition

\remark{Remark 18} Let a family $\fe$ be nondegenerate. Then the multipliers 
$(\fe)'(\aje)$ of the maps $\fe$ at their fixed points $\aje$ tend to 1, as 
$\var\to0$, and the differences $\fe'(\aje)-1$ have arguments bounded away from 
$\frac{\pi}2+\pi\Bbb Z$ uniformly in  all $\var$ 
small enough. In particular, for all $\var$ small enough the multipliers 
lie in the complement of unit circle. 
\endremark

\define\ve{v_{\var}}

\proclaim{Theorem 7} Let $\fe$ be a nondegenerate family of conformal maps, 
$\a=\a_i(\var)$ be a continuous family of their fixed points, $S=S_{j_i}$ be a
correspondent sector from either the item following Remark 8,  
if $k=1$ (then $S_{j_i}=S_i$), or Remark 11 otherwise. There exist an $r>0$ and a family of domains  
$\ov=\ov(i)$  depending on the same parameter containing $\a(\var)$ and no other fixed 
point with the following properties:  

 1) The connected component of the intersection $\ov\cap (S^r\setminus\bigcup_{\a_s\neq\a}[0,\a_s(\var)])$ containing $\a(\var)$ converges to $S^r$, as $\var\to0$.
 
2) The family $\ve$ of the canonic generators 
of the perturbed maps at $\a(\var)$ is well-defined  for all 
parameter values small enough and each $\ve$ extends analytically to $\ov$;
 the family $\ve$ converges to the sectorial canonic generator of 
the nonperturbed map in $S^r$, as $\var\to0$.
\endproclaim

Theorem 7 is proved in Subsection 3.4.

\proclaim{Corollary 3} Let $\fe$ be a nondegenerate family of conformal maps, 
$\a_i$ be its continuous fixed point families, $S_{j_i}^r$, $\ov(i)$ be  respectively the correspondent sectors and 
domains from Theorem 7. Let $\tau_{i,\var}(t)$ be 
the complex time functions correspondent to the canonic generators of $\fe$ at $\a_i(\var)$. 
The functions $\tau_{i,\var}$ are holomorphic and single-valued in the complements $\ov'(i)=\ov(i)
\setminus\cup_{s=0}^{k}[0,\a_s(\var)]$. The families $\tau_{i,\var}$ of 
appropriately normalized time functions converge to the sectorial time functions   $\tau_{j_i}$ from Remark 16 in $S_{j_i}^r$. Let  $\a_i$ and $\a_{i+1}$ be neighbor fixed point families such that the 
correspondent sectors intersect each other, i.e., $S_{j_{i+1}}=S_{j_i+1}$ (we put $\tau_{2k+1}=\tau_1+\lambda$, $\a_{k+1}=\a_0$, $\tau_{k+1,\var}=\tau_{0,\var}+\lambda$).    There exists a family $C_\var$ of connected components of the intersection  
$S_{j_i}^r\cap S_{j_{i+1}}^r\cap\ov'(i)\cap\ov'(i+1)$ that tends to the  component of the intersection $S_{j_i}^r\cap S_{j_{i+1}}^r$, as $\var\to0$. Let $\tau_{i+1, \var}\circ\tau_{i,\var}^{-1}$ be the  transition function  
that compares the correspondent time functions in $C_\var$. The transition function is well-defined in a domain (depending on $\var$) that tends to a half-plane $(-1)^{j_i}\operatorname{Im}\tau<c$. The family of the transition functions converges to the  correspondent Ecalle - Voronin transition function $\psi_{j_i}$ of the nonperturbed map $f_0$ (see Remark 16, we put $\psi_{2k}=\psi_0$). \endproclaim

A particular case of Theorem 7 for $k=1$ and nondegenerate 
deformations $\fe$ analytic in $\var$ such that $(\fe)'_{\var}(0,0)\neq0$ was proved in 
\cite{4}.

\head 3. Proof of Theorems 3, 4, 7 and Lemma 1\endhead

\subhead 3.1. Scheme of the proof of Theorems 3 and 4 \endsubhead

We prove Theorems 3 and 4 in Subsections 3.1-3.3 and 3.5. Another their proof that uses Theorem 7 and is more short is sketched 
in Subsection 3.6. 

 For the proof of Theorems 3 and 4 it suffices to show that there exists a 
 family $\tz=\he(z,t)$ of changes of the variable $z$ over neighborhoods $\ov$ 
 of $\ae$ satisfying statement 1) of Theorem 3 (4) that linearizes the differential 
 equation in $z$ of $\le$, i.e., transforms it to a family of equations of 
 the type 
 $$\dot{\tz}=\tz g_{\var}(t),$$
and converges to the sectorial normalization $H(z,t)$ from Theorem 1 of the nonperturbed field. Then the canonic first integral 
of the perturbed field takes the form $$\tz e^{-\int^t\frac{g_{\var}(\tau)d\tau}{\prod_j(\tau-\aje)}},$$
and by definition, its restriction to the correspondent domain $\widetilde{\ov'}$ 
from Theorem 3(4) will converge to the sectorial canonic first integral of the 
nonperturbed field. 

\define\go{\Gamma_0}

\definition{Definition} In the conditions of Theorem 3(4) 
the {\it sectorial separatrix} $\go$ of the nonperturbed field over the sector $S$ is the zero level curve of the correspondent sectorial canonic integral.  
 The separatrix $\ge$ of the perturbed field is the 
zero level curve of its canonic integral at $(0,\a(\var))$. 
\enddefinition

\remark{Remark 19} Let (1) be a vector field as in Remark 1, $U_z$ be a  neighborhood of zero in the $z$- axis. There exists an $r>0$ such that the sectorial separatrix $\go$ of the nonperturbed field contains the graph $z=q(t)$, $t\in S^r$, of a $U_z$- valued function $q$ holomorphic in the sector $S^r$ and continuous in its closure, $q(0)=0$. This is the unique phase curve of the nonperturbed field that contains such a graph. 
This follows from Theorem 1 and the fact that this is valid for the formal normal form. 
\endremark

\remark{Remark 20} Let $\le$ be a nondegenerate vector field family, $(0,\a(\var))$ be its continuous singularity family, $U_z$ 
be a fixed neighborhood of zero in the $z$- line. 
For any $\var\neq0$ the separatrix $\ge$ of the perturbed field at the singular point $(0,\a(\var))$ extends analytically to the latter and is tangent to  the eigenline of the correspondent linearization operator transversal to the line $t=\a(\var)$ (the correspondent eigenvalue tends to 0, as $\var\to0$). In particular, it contains 
the graph $z=\qe(t)$ of a function $\qe$ holomorphic in a neighborhood of 
$\a(\var)$ (that depends on $\var$), $\qe(\a(\var))=0$. 
\endremark

For the proof of the statement from the beginning of the Subsection on the existence and convergence of a linearization family $\he$ we firstly show that 
the separatrices $\ge$ of the perturbed fields converge to the sectorial 
separatrix $\go$ of the nonperturbed field: 

\define\dtd{\{|t|<\dl\}}
\proclaim{Lemma 2} Let $\le$, $\ae$, $S$ be as in Theorem 3 (4). There exist 
an $r>0$ and neighborhoods $U_t=\dtd$, $U_z=\{|z|<2\dl\}$ of zero 
in the $t$- and $z$- lines respectively such that there exists a family $\ov^1$ 
of subdomains in $U_t$ containing $\ae$ and no other $\a_j(\var)$ 
satisfying statement 1) of Theorem 3 (4) with the following property:  
 for all $\var\neq0$ small enough  the separatrix $\ge$ from the previous Definition contains the graph of a holomorphic function 
$z=\qe(t)$ in $t\in\ov^1$ with values in $U_z$. The family $\qe$ converges to the function $q$ from Remark 19 correspondent to the sectorial separatrix $\go$ of the nonperturbed field.  
\endproclaim 

 Lemma 2 is proved in Subsections 3.2 and 3.5. (A proof of its version is also implicitly contained in \cite{6}.) The correspondent domains $\ov^1$ are defined in the next Subsection (after Remark 21). In the case, when, $k=1$, they are  depicted at Fig.3a. They converge to a cardioid-like domain bounded by a Jordan curve with cusp at 0 depicted at Fig.3b.

In the proof of the statement from the beginning of the Subsection we use Lemma 2 and the following statement on the "uniqueness" of the sectorial $z$- variable  linearizing chart for the nonperturbed field.

\proclaim{Proposition 1} In the conditions of Theorem 1 let 
$U=\{|z|<2\delta\}\times\{|t|<\delta\}$, $S=S_j$ be a good 
sector. Let $\g_0$ be the zero level curve of the correspondent sectorial canonic integral. Let $H(z,t)$ be a holomorphic function in $\ts$ univalent in $z$ in the 
discs $t=const$ and vanishing in $\g_0$ such that  the change $\tz=H(z,t)$ of 
the variable $z$ linearizes the differential equation in $z$ of (1). Then $H$ 
 is obtained from the correspondent sectorial normalizing function 
$H_j$ by multiplication by holomorphic function in $t$.    
\endproclaim

\demo{Proof} It suffices to prove Proposition 1 for the formal normal form $(1)_n$ (Theorem 1), where the equation in $z$ is already linear. 
Let $H$ be as in Proposition 1. Let us show that $H$ is 
linear in $z$. This will prove Proposition 1. 

By definition, $H(0,t)\equiv0$. Without loss of generality we consider that 
$\frac{\partial H}{\partial z}(0,t)\equiv1$. One can achieve this by  
multiplying $H$ by appropriate function in $t$. 

We prove Proposition 1 by contradiction. Suppose the function $H$ is not 
linear in $z$. Then it takes the form 
$$H(z,t)=z(1+z^lg(t)+\text{higher terms in}\ \ z).$$
The function $g$ is bounded, which follows from univalency in $z$ of the function 
$H$ and distorsion theorem from \cite{10}. On the other hand, it satisfies the 
equation 
$$\frac{dg}{dt}=\frac l{t^{k+1}}g,$$
which follows from the conditions of 
Proposition 1 (the differential equation of the field 
in the new variable $\tz=H(z,t)$ is linear, as that in $z$). Therefore, $g$   is a constant multiple of the function $e^{\frac l{kt^k}}$,
which tends to infinity exponentially, as $t\to0$ 
along a ray where the real part of its exponent is positive. The sector $S$, which is good, contains such a ray by definition: it contains an  
imaginary dividing ray, where this real part changes sign. This contradicts the boundedness of $g$. Proposition 1 is proved.
\enddemo

As it is shown in Subsections 3.3 nd 3.5 (Lemma 3 from the end of the present Subsection), there exist a disc $U_z$ in the $z$- axis and $z$- variable linearizing charts $\tilde z=\he(z,t)$ for the perturbed fields (which vanish at $\ge$) over appropriate domains $\ov^2$ in the $t$- line (containing $\a(\var)$ and satisfying statement 1) of Theorem 3 (4)) univalent in $z\in U_z$ for any fixed $t\in\ov^2$. At the end of Subsection 3.3 we prove that the connected component $\ov$ of the intersection $\ov^1\cap\ov^2$ ($\ov^1$ is the domain  from Lemma 2) that contains $\a(\var)$ satisfies statement 1) of Theorem 3 (4). This is the domain family we are looking for. One can choose such a family $H_\var$ of linearizing charts 
to have unit derivatives in $z$ at $\ge$. Then the family $\he|_{\widetilde{\ov}\cap \ts^r}$ is normal with respect to the compact convergence in $\widetilde{S^r}$ (Lemma 2 and normality of the space of the normalized univalent functions  \cite{10}). Using  normality and Proposition 1, let us show that this family converges to the sectorial normalizing chart of the nonperturbed vector field over $S^r$ (multiplied by a nonvanishing holomorphic function in $t$). This will prove Theorems 3 and 4.  

By construction, the limit of each convergent sequence $H_{\var_n}$, $\var_n\to0$, is a $z$- linearizing chart over $S^r$ for the nonperturbed field  that satisfies the conditions of Proposition 1 and has unit derivative in $z$ at $\go$. By Proposition 1, this chart is obtained from the canonic normalizing chart by multiplication by nonvanishing holomorphic function in $t$. The chart 
satisfying conditions of Proposition 1 and having unit derivative in $z$ at $\go$ is unique, which follows immediately from Proposition 1. In particular, all the limits of convergent sequences coincide, and hence, the family $\he$ 
converges to a chart as at the end of the previous item. This proves Theorems 3 and 4 modulo Lemma 2, the following Lemma and 
statement 1) of Theorem 3 (4) for the correspondent domains $\ov$ from the previous item. The last statement is proved at the end of Subsection 3.3.

\proclaim{Lemma 3} In the conditions of Lemma 2 there exist an $r>0$,  
neighborhoods $U_z$, $U_t$ of zero in the $z$- and $t$- axes respectively, a  family of subdomains $\ov^2\subset U_t$ satisfying statement 1) of Theorem 3 (4)  such that there exists a family $\he(z,t)$ of functions holomorphic in $\widetilde{\ov^2}=U_z\times\ov^2$
and univalent in $z\in U_z$ in the discs $t=const$ such that  the coordinate change 
$\tz=\he(z,t)$ linearizes the differential equation in $z$ of the perturbed 
field. 
\endproclaim

Lemma  3 is proved in Subsections 3.3 and 3.5. 

\subhead 3.2. Continuity of the separatrices. Proof of Lemma 2\endsubhead

\define\vt{v_\theta}

\define\wt{w_\theta}

In the proof of Lemma 2 we use the uniqueness statement of Remark 19 from the 
previous Subsection.  

Let $q$ and $\qe$ be the functions from Remarks 19 and 20 correspondent respectively to the sectorial separatrix 
of the nonperturbed vector field and the separatrix of the perturbed field. 
For the proof of convergence $\qe\to q$ in appropriate sector $S^r$ we show that the function family $\qe$ is equicontinuous in a family of domains convergent to $S^r$. Then they form a normal family (Arzela-Ascoli theorem).  
Each limit of sequence $\{ q_{\var_n}\}$, $\var_n\to0$,  compactly convergent
in $S^r$ is holomorphic in $S^r$ and extends continuously to its closure. The graph of each limit is contained in a phase curve of the limit (nonperturbed) field. Therefore, by the uniqueness statement of Remark 19, these limits coincide with 
the function $q$ correspondent to the sectorial separatrix $\go$. This will 
prove Lemma 2.    

The equicontinuity statement from the last item (and hence, Lemma 2 also) is implied by the following 

\proclaim{Lemma 4} In the conditions of Lemma 2 there exist an $r>0$, 
a neighborhood $U_z$ of zero in the $z$- line and a family $\ov^1$ of domains 
in the $t$- line containing $\a(\var)$ and satisfying statement 1) of Theorem 3 (4) such that the separatrix $\ge$ of the perturbed vector field 
at $(0,\a(\var))$ contains the graph $z=\qe(t)$ of a function $\qe(t)$, $\qe(\a(\var))=0$, holomorphic 
in $\ov^1$ with values in $U_z$ that satisfies the following inequalities:
$$|\frac{d\qe}{dt}(t)|<1,\ \ \ |\qe(t)||_{\a(\var)\neq t\in\ov}<|t-\ae|.\tag3$$
\endproclaim

\demo{Proof} 
Before the proof of Lemma 4, let us introduce the following 

\definition{Definition} Let (1)$_\var$ be a nondegenerate family,  
$\theta\in\Bbb R$. Define $\vt(\var)$ to be the vector field family that is the 
$e^{i\theta}$-th multiple of (1)$_\var$. Define $\wt(\var)$ to be the projection  image in the $t$- line of the restriction $\vt(\var)|_{z=0}$. 
\enddefinition

 The inequalities of Lemma 4 are equivalent to the statements 
that the tangent lines to $\ge$ are contained in the tangent cone field 
 $K=(|\dot t|>|\dot z|)$ and $\ge$ is contained in the cone $\widetilde K=
 \{|t-\ae|>|z|\}$: 
 
$$T\ge\subset K,\ \ \ \ \ \ \ \ \ \ge\subset\widetilde K.\tag4$$  

Inclusions (4) hold a priori in a neighborhood of the singular point $(0,\alpha(\var))$ 
(that depends on $\var$). To show that they hold in a large domain, 
we consider a vector field family $\vt(\var)$ from the previous Definition with the following properties:

\item{1)} the eigenvalue of the linearization operator of the field $\vt(\var)$ at the eigenline tangent to the line $t=\ae$ has negative real part bounded away  from 0 uniformly in $\var$ small enough; 

\item{2)} the other eigenvalue has positive real part, and its argument is bounded away from $\frac{\pi}2+\pi\Bbb Z$ uniformly in $\var$ small enough. 

\remark{Remark 21} Let (1)$_\var$ be a nondegenerate family, $\alpha(\var)$ be 
its continuous singularity coordinate family. There exists a 
$\theta\in\Bbb R$, $\frac{\pi}2<\theta<\frac{3\pi}2$, that satisfies the conditions 
1) and 2) from the two previous items. Indeed, for $\theta$ satisfying the last 
inequality condition 1) holds by definition. Condition 2) is satisfied, e.g., 
for $\theta_0=\pi\pm\frac{\pi}2$ with appropriately chosen sign 
(dependently on the choice of the family $\ae$), which follows from nondegeneracy of (1)$_\var$. Thus, 2) holds also for all $\theta$ close enough to $\theta_0$, which may be chosen to satisfy the previous inequality. Under 
the conditions 1) and 2) $(0,\ae)$ is a hyperbolic singularity of the field 
$\vt(\var)$, with 1-dimensional stable and unstable manifolds: the former is 
contained in the line $t=\ae$; the latter is contained in $\ge$. The singular 
point $\ae$ of the field $\wt(\var)$ is repelling. 
\endremark

We use the fact that under conditions 1) and 2) there exists a neighborhood 
$U$ of zero in the phase space such that for all $\var$ small enough the cone 
field $K$ (and hence, $\widetilde K$ as well) is $\vt(\var)$- invariant in $U$ 
(proposition 2 in \cite{6}). This together with the local inclusions (4) in a 
neighborhood of the singular point $(0,\ae)$ implies that the latters hold in the trajectories of $\vt(\var)$ in 
 $\ge\cap U$ that go from $(0,\ae)$ (denote by $\ov''$ the subdomain of the separatrix  
 $\ge$ saturated by these trajectories). Without loss of generality we consider that $U=U_z\times U_t$, $U_z=\{|z|<2\delta\}$, $U_t=\{|t|<\delta\}$. Then inclusions (4) hold over the following domains.

\definition{Definition of the domains $\ov^1$} In the conditions of Theorem 3(4) and the Definition from the beginning of the proof of Lemma 4 let $\theta$ satisfy the condition 2) preceding Remark 21. Let $U_t$ be a fixed disc centered at  zero in the $t$- line. Define $\ov^1\subset U_t$ to be the subdomain saturated by the trajectories of the 
 field $\wt(\var)$ in $U_t$ that go from $\ae$. (In the case, when $k=1$, the phase portrait of the field  $\wt(\var)$ 
 and the domain $\ov^1$ are depicted at Fig.3a. The boundary of the domain $\ov^1$ consists of the three following parts: a) the arc of the boundary $\partial U_t$ where the field $\wt(\var)$ is directed outside $U_t$ bounded by 
the pair of tangency points of the field with $\partial U_t$; b) the pair of semitrajectories of the field that start at these points and tend to its other 
singularity.) 
\enddefinition

Inclusions (4) hold 
 in $\ov''$ (see the item preceding the Definition). The latter 
 is 1-to-1 projected onto $\ov^1$. This follows from the inclusion 
 $\ov''\subset\widetilde K$ and the definition of the neighborhood $U$ (exactly in the same way, as in the proof of lemma 1 in \cite{6}). Therefore, $\ov''$ is the graph of a function $q_\var(t)$ holomorphic in $\ov^1$ and satisfying inequalities (3). 
 
 Now for the proof of Lemma 4 it suffices to show that one can choose $\theta$ satisfying conditions 1) and 2) so that 
 the correspondent family $\ov^1$ satisfies statement 1) of Theorem 3 (4). To do this, 
we prove that appropriate family of subdomains of $\ov^1$ converges to a domain $\Omega^1$ bounded by a 
Jordan curve having cusp at 0 of the following type. In the case, when $k=1$, the branches of the cusp are tangent to 
the ray $-R(\theta)$ with the argument $\pi-\theta$ disjoint from $\Omega^1$  (Fig.3b). We show that one can choose $\theta$ satisfying conditions 1) and 2) so that the correspondent domain $\Omega$ contains the sector $S^r$ for appropriate $r>0$. To do this, we prove the possibility to choose $\theta$ as above so that the ray $-R(\theta)$ punctured at 0 lies outside the closure of the sector $S$. In the case, when, $k\geq2$, the cusp branches of $\partial\Omega^1$ are tangent to a pair of rays that bound a sector disjoint from $\Omega^1$. We show that one can choose $\theta$ satisfying conditions 1) and 2)  so that the closure punctured at 0 of the last sector is disjoint from that of $S$. Then $\Omega$ will contain $S^r$ for appropriate $r>0$. 
  
 \definition{Definition} Let $k\in\Bbb N$, $k\geq1$, $\theta\in\Bbb R$. Define $R_m(\theta)$ to be the radial ray with the argument $\frac{-\theta+\pi m}k$.
 \enddefinition

In the proof of the statements from the item preceding the previous Definition 
 we use the following properties of the nonperturbed vector field $\wt(0)$. 
 
 \remark{Remark 22} In the conditions of the previous Definition each trajectory of the vector field $\dot t=e^{i\theta}t^{k+1}$ has one of the three 
 following types: (i) a ray $R_m(\theta)$ with $m$ even oriented from 0 to 
 infinity (repelling ray); (ii) a ray $R_m(\theta)$ 
 with $m$ odd oriented towards 0 (attracting ray); (iii) a Jordan curve that 
 goes from zero with asymptotic tangency to a ray from (i) and returns to 0 with  the asymptotic tangency to a neighbor ray from (ii). \endremark 
 
 \proclaim{Corollary 4} In the conditions of the previous Remark let $R_m$ be a repelling ray, $U_t$ be a 
 disc centered at 0. Then the subdomain $\Omega_m\subset U_t$ saturated by the trajectories of the field in $U_t$ that go from 0 with asymptotic tangency to the ray $R_m$ 
 is bounded by a Jordan curve that has a cusp at 0 with the two tangency rays 
 $R_{m\pm1}$ disjoint from $\Omega_m$. (This subdomain is called the repelling domain of $R_m$. In the case, when $k=1$, it is depicted at Fig.3b.) Its boundary consists of the three 
 following parts: an arc of $\partial U_t$ where the field is directed outside 
 $U_t$ bounded by tangency points of the field with $\partial U_t$; two  semitrajectories that start at these points and converge to 0 with tangency 
 to the rays $R_{m\pm1}$. 
\endproclaim 

\example{Example 3} Let $k=1$. Then in the conditions of the previous Remark the  only repelling ray $R(\theta)$ is that with the argument $-\theta$, and the 
opposite one ($-R(\theta)$) is 
the only attracting ray. The correspondent repelling domain $\Omega_m$ from Corollary 4 is 
bounded by a cardioid-like Jordan curve having inward cusp at 0 with 
tangency to the attracting ray, which is disjoint from $\Omega_m$ (Fig.3b).
\endexample

\hskip-2cm Fig.3 a,b

 For appropriate repelling domain $\Omega^1$ from 
Corollary 4 correspondent to the nonperturbed field $\wt(0)$ we show that the domains $\ov^1\cap\Omega^1$ converge to $\Omega^1$. This is implied by the following 

\proclaim{Proposition 2} Let $\theta\in\Bbb R$. Let 
$w(\var):\ \dot t=e^{i\theta}(t-\a_0(\var))\dots(t-\a_k(\var))$ be a continuous 
family of vector fields 
in complex line with the coordinate $t$ having $k+1$ distinct singularities $\a_i(\var)$ 
that tend to 0, as $\var\to0$. In the case, when $k\geq2$, let the singularity polygon satisfy the asymptotic regularity  statement of Remark 10. Let $U_t$ be a neighborhood of zero in the phase line.  Let 
$\ae=\a_j(\var)$   be a continuous singularity family. Let
 the multiplier of the field $w(\var)$ at $\ae$ have positive real part and its 
argument be bounded away from $\frac{\pi}2+\pi\Bbb Z$. Let $\ov\subset U_t$ be the correspondent family of subdomains from the Definition following Remark 21.  Then 
there exists a repelling ray $R$ (see the previous Remark) correspondent to the 
nonperturbed field $w(0)$ such that for all $\var$ small enough $R$ is the 
closest of all the repelling and attracting rays to the radial ray of $\ae$. 
Let $\Omega$ be the correspondent repelling domain from Corollary 4. In the case, when $k=1$, the family $\ov$ tends to $\Omega$, as $\var\to0$. In the 
case, when $k\geq2$, so does the family of the connected components of the intersections $\ov\cap\Omega$ containing $\a(\var)$. In both cases there exists 
an $r'>0$ such that each $\ov$ contains the segment of the radial ray of the point  $\a(\var)$  joining the latter to the point of the circle $|t|=r'$. 
\endproclaim

A proof of Proposition 2 except for that of its two last statements is implicitly contained in subsections 5.A and 6 of the paper \cite{6}. Its last statements are proved  in Subsection 3.5. 

For simplicity, we consider only the case, when $k=1$, and prove the existence of $\theta$ such that $\ov^1$ satisfy 
statement 1) of Theorem 3. The proof remains 
valid for larger $k$ with obvious changes. 

 We show that one can choose $\theta$ satisfying conditions 1) and 2) from 
the item preceding Remark 21 so that the attracting ray $-R(\theta)$ from the previous Example punctured at 0 
lies outside the  closure of the sector $S$ correspondent 
to $\ae$. Then the correspondent "limit repelling domain" from 
Proposition 2 (it will be referred to, as $\Omega^1$) will contain a sector 
$S^r$ for some $r>0$. Thus, a connected component of the intersection from statement 1) of Theorem 3 tends to $S^r$. This component contains $\a(\var)$. 
Indeed, for all $\var$ small enough it contains the segment from the last statement of Proposition 2 correspondent to $r'<r$. This follows from the fact that the segment end different from $\a(\var)$ is contained in the sector $S^r$ and bounded away from its boundary (by definition) uniformly in all small $\var$. This proves statement 1) of Theorem 3 modulo the two last statements of Proposition 2. 

It suffices to prove the first statement from the last item for the 
argument $\theta_0$ from Remark 21 instead of $\theta$, since $\theta$ satisfying conditions 1) and 2) may be chosen arbitrarily close to $\theta_0$ 
(Remark 21). The statement for $\theta_0$ follows 
from the fact that $R(\theta_0)$ is the imaginary dividing ray that lies in 
$S$ (so, the opposite imaginary dividing ray $-R(\theta_0)$ lies outside the sector $S$, since the latter is good). This follows from the definitions of $\theta_0$ (Remark 21) and the correspondent ray $R(\theta_0)$ (Example 3). Statement 1) of Theorem 3 together with Lemma 2 is proved modulo the two last statements of Proposition 2.  
\enddemo

\subhead 3.3. Proof of Lemma 3. The end of the proof of Theorems 3 and 4 modulo Proposition 2\endsubhead 

Firstly let us prove Lemma 3. In its proof we use the following 

\remark{Remark 23} Let 
$$\cases
& \dot z=z+O(|z|^2+|t|^2) \\
& \dot t=\lambda t(1+g(t)),\endcases$$   be a holomorphic vector field in a 
neighborhood of zero in $\Bbb C^2$, $g(0)=0$, $\lambda\notin\Bbb R$. There exists a unique change $\tilde z=H(z,t)$ of the variable $z$ in a neighborhood of zero that linearizes in $z$ the correspondent (first) differential equation  up to multiplication by nonzero function in $t$. 
\endremark

 A $z$- linearizing chart $\tilde z=\he(z,t)$ for the perturbed vector field exists a priori in a neighborhood $U_z'\times U_t'$ of 
the singular point $(0,\ae)$ (that depends on $\var$). Without loss 
of generality we consider that it is univalent in $z$ in the discs $U_z'\times t$ and has unit derivative in $z$ at $\ge$. To show that it continues to a large domain, we consider a vector 
 field $\vt(\var)$ (correspondent to another value $\theta$) that satisfies 
 the conditions 1) and 2) preceding Remark 21 from the previous Subsection with the 
 change of "negative real part" in the condition 1) to "positive real part". 
 Then its singular point $(0,\ae)$ is repelling, i.e., it is asymptotically 
 stable with respect to the opposite field $-\vt(\var)$. Denote by $g^{\tau}$ the time $\tau$ 
 flow map of the field $\vt(\var)$. We use the fact that for any $\tau$ and 
 linearizing function $\he$ (defined in a neighborhood of the singularity) 
 the composition $H_{\var,\tau}=\he\circ g^{\tau}$ also defines a change of the coordinate $z$ 
 linearizing the correspondent differential equation. Thus, it is obtained from 
 $\he$ by multiplication by function in $t$ (the previous Remark). Therefore, by definition, one has 
 $\he=\frac{H_{\var,\tau}}{((H_{\var,\tau})'_z)|_{\ge}}$. Let $U=U_z\times U_t$ 
 be a fixed bidisc centered at zero in the phase space. The last formula defines the analytic 
 continuation of $\he$ along the trajectories of $\vt(\var)$ in $U_t$ that go from $(0,\ae)$, thus, to the domain saturated 
 by the latters. This domain will be referred to, as $\ov''$. For appropriate neighborhoods $U_z$ and $U_t$ and all $\var$ small enough $\ov''$ is 
the Cartesian product $\widetilde{\ov^2}$ of the disc $U_z$ and the domain from the previous Definition (it will be referred to, as $\ov^2$) correspondent to the new $\theta$: namely, if the field $-v_{\theta}(\var)$ is directed inside the domain $\widetilde{\ov^2}$. This is the case, if the field $-v(\theta)$ is directed inside the 
 cylinder $U_z\times \Bbb C$ over $U_t$ (by definition). The possibility of choice of neighborhoods $U_z$ and $U_t$ satisfying the last condition for all $\var$ small enough 
 follows immediately from the choice of $\theta$ (the modified condition 1)) and the continuity of the family $\le$. Then by construction, the analytic continuation of the function $\he$ to $\widetilde{\ov^2}$ is linearizing and univalent in $z$. This  
follows from definition and the fact that the 
 fibration by parallel lines $t=\operatorname{const}$ is preserved under the 
 flow $g^\tau$. This proves Lemma 3 for the domains $\ov^2$ modulo statement 1) of Theorem 3 (4). 
 
 Similarly, as in the proof of Lemma 2, one can choose $\theta$ as in 
 the previous Definition so that the family $\ov^2$ satisfies statement 1) of Theorem 3 (4), 
 more precisely, the correspondent repelling domain from Proposition 2 contains the sector  $S^r$ for appropriate $r>0$. This proves Lemma 3 (modulo the last statements of Proposition 2).  
 
 Now Theorems 3 and 4 are proved modulo statement 1) of Theorem 3 (4) for the following domains $\ov$  (see the discussion in Subsection 3.1 preceding the statement of Lemma 3).
   
 \definition{Definition of the domains $\ov$} In the conditions of Theorem 3(4) and the Definition following Lemma 4 let $\theta_1, \theta_2\in\Bbb R$ satisfy condition 2) preceding Remark 21, $\theta_1$ satisfies condition 1) from the same place, 
$\theta_2$ satisfies its modification with the change of "positive real part" to  "negative real part". Let $U_t$ be a disc centered at 0 in the $t$- axis. Let $\ov^1$, $\ov^2$ be the domains from the Definition of the domains $\ov^1$ in the previous Subsection correspondent to $\theta_1$ and $\theta_2$ respectively.  Define $\ov$ to be the connected component of the 
 intersection $\ov^1\cap\ov^2$ that contains $\ae$. (In the case, when $k=1$, the boundary of the domain 
 $\ov$ consists of the three following parts: an arc of $\partial U_t$ where 
 both fields $w_{\theta_1}(\var)$ and $w_{\theta_2}(\var)$ are directed 
 outside $U_t$ bounded by the tangency points $F_1(\var)$ and $F_2(\var)$ 
 respectively of these fields with $\partial U_t$; 
 two arcs of their semitrajectories (tangent to $\partial U_t$) that  start at $F_1(\var)$ and $F_2(\var)$ respectively and have a common end, see Fig.4a.) 
 \enddefinition 

\hskip-2cm Fig. 4 a,b

For appropriate values $\theta_1$ and 
$\theta_2$ from the previous Definition the correspondent family $\ov$  satisfies statement 1) of Theorem 3 (4). This is the case, provided that 
so do $\ov^1$ and $\ov^2$ (more precisely, the correspondent limit domains $\Omega^1$ and $\Omega^2$ from Proposition 2 contain the sector $S^r$ for appropriate $r>0$). Indeed let $\Omega$ be the connected component of the intersection 
 $\Omega^1\cap\Omega^2$ that contains the sector $S^r$. The domain $\Omega$ is bounded by a Jordan curve that consists of the three following parts: an arc of 
$\partial U_t$ where both $w_{\theta_1}(0)$ and $w_{\theta_2}(0)$ are directed 
outside $U_t$ bounded by the tangency points $F_1$ and $F_2$ respectively 
of the latters with $\partial U$; 
a pair of positive semitrajectories $L_1$ and $L_2$ of these fields that start at the 
points $F_1$ and $F_2$ respectively and converge to 0 (Fig.4b). In the case, 
when $k=1$,   
$\ov\to\Omega$, as $\var\to0$, which was implicitly proved in subsection 5.B of \cite{6}. In the case, when $k\geq2$,  the connected component of the intersection $\ov\cap\Omega$ containing $\a(\var)$ tends to $\Omega$, as $\var\to0$. Indeed, by the convergence statement of Proposition 2 applied to  $\ov^i\cap\Omega^i$, $i=1,2$, some connected component of the intersection 
$\ov\cap\Omega$ converges to $\Omega$. This together with the last statement of Proposition 2 implies that this component contains $\a(\var)$, whenever $\var$ small enough, as in the discussion preceding the last item of Subsection 3.2. Analogously statement 1) of Theorem 3 (4) follows from the convergence statements for the domains $\ov$ and the last statement of Proposition 2 for the domains $\ov$ (it holds for both $\ov^i$). Theorems 3 and 4 are  proved modulo the two last statements of Proposition 2.

\define\mtd{\{|t|<\delta\}}

\subhead 3.4. Proof of Theorem 7 \endsubhead

The method of the proof of Theorem 7 is similar to that of Theorems 3 and 4. 

Without loss of generality we consider that in the conditions of Theorem 7 
$\ae$ is a repelling fixed point family for $\fe$. One can achieve this by 
changing $\fe$ to their inverse maps and subsequent linear coordinate change. The generator $\ve$ of the perturbed map $\fe$ is a priori defined in a 
neighborhood of $\ae$ (that depends on $\var$) and is $\fe$- invariant. It  extends analytically to 
the domain $\ov$ from the next Definition. The continuation is made along the orbits of $\fe$  that go from 
$\ae$ by applying the iterations of $\fe$ to the generator defined near $\ae$. 

\definition{Definition of the domains $\ov$} Let $\fe$ and $\ae$ be as in the last item and 
Theorem 7. Let $U=\mtd$ be a fixed neighborhood of zero where $\fe$ are 
holomorphic and have the only fixed points $\a_j(\var)$ for all $\var$. 
For any $\var$ define $\ov\subset U$ to be the subdomain saturated by the 
orbits of the map $\fe$ in $U$ that go from $\ae$ (see the reversed 
Fig.3a for $k=1$). 
\enddefinition

We show that the domains $\ov$ satisfy the statements of Theorem 7. 

 Without loss of generality, we consider that the argument family of the fixed points $\a(\var)$ has a limit, as $\var\to0$ (in the case, when $k\geq2$, this follows from the definition of nondegeneracy). Each parameter value sequence 
$\var_n\to0$ contains a subsequence where the arguments $\arg\a(\var)$ converge. Therefore, the convergence of the generators $\ve$ of the perturbed maps in these subsequences will imply their convergence in the whole parameter space. 

By $p(t,\var)$ denote the vector field family $\dot t=\pr$. 

Firstly we prove statements of Theorem 7 for a smaller good sector $C\subset S$.  For the proof of convergence of the generators $\ve$ we show that there exist $r,r',b>0$, a family of $\fe^{-1}$- invariant subdomains $D_\var\subset\ov$ containing $\a(\var)$, and a good sector $C\subset S$ with the following properties: 

\item{1)} the family $D_\var$ converges to a (connected) domain $D_0$ containing  the sector $C^r$; 

\item{2)} the inequality $|\ve(t)|\leq b|p(t,\var)|$ holds in $D_\var$ for all $\var$ small enough;

\item{3)} for any $\var$ small enough the domain $D_\var$ contains the segment 
of the radial ray of the point $\a(\var)$ joining the latter to the point of 
the circle $|t|=r'$.  

(In the case, when $k=1$, for $\var\neq0$ the domain 
$D_\var$ will be a disc symmetric with respect to the line passing through the 
fixed point pair of the map $\fe$, the limit domain $D_0$ will be the maximal  disc in $U$ symmetric with respect to this line whose boundary contains 0 (Fig. 5a,b).)  
Then the family $\ve$ converges to the sectorial generator $v_0$ of the nonperturbed map in $D_0$ (in particular, in $C^r$). This  is proved by  using the bound 2), normality of the family $\ve$ (implied by 2)) and uniqueness of the sectorial generator (Remark 15), analogously to the proof of Lemma 2 (the discussion preceding Lemma 4 in Subsection 3.2). 

\redefine\we{w_{\var}}

Now let us show how the convergence statement in the smaller good sector $C$ together with the previous statement 3) implies Theorem 7 (for the larger sector $S$). Let $r_j$ be the imaginary dividing ray contained in the sector $S$, which has the argument $\frac{\pi}{2k}+\frac{\pi j}k$. Consider the subdomain $\Omega\subset U$ saturated by the orbits of the limit map $f_0$ in $U$  that go from zero in the asymptotic direction of the ray $r_j$. The domain $\Omega\subset U$ is bounded by Jordan curve that has cusp at 0 with tangency to the two imaginary dividing rays neighbor to $r_j$. This is  implicitly proved in \cite{1}. (In the case, when $k=1$, these two rays coincide, see the  reversed Fig.3b.) In particular, $\Omega\supset S^r$ for appropriate $r>0$. A connected component of the intersection $\ov\cap\Omega$ converges to $\Omega$, as $\var\to0$. This follows from the definition of the domains $\ov$, the property 1) of the domains $D_\var$ and the fact that $\Omega$ is saturated by the $f_0$- orbits that start in $C^r$ (by definition). (In the case, when $k=1$, $\ov\to\Omega$. The proof of this statement is omitted to save the space.) Therefore, the same convergence statement holds with the change of $\ov\cap\Omega$ to the intersection from statement 1) of Theorem 7. The sectorial generator $v_0$ of the nonperturbed map extends analytically from $C^r$ to $\Omega$ by applying iterations of the map $f_0$ to its restriction to $C^r$. The analogous continuation of the generator $\ve$ to $\ov$ along the $\fe$- orbits starting in $D_\var\cap C^r$ converges to $v_0$ in $\Omega$. This follows from continuity of the family $\fe$ and the convergence of $\ve$ in $C^r$. This proves Theorem 7 for the sector $S^r$ modulo the statement that the convergent component of the correspondent intersection contains $\a(\var)$. This follows from the convergence of the component under consideration and the previous property 3) of $D_\var$ (analogously to the discussion preceding the last item of Subsection 3.2). 

In the proof of the previous statements 1)-3) for appropriate domains $D_\var$ we use the following 

\remark{Remark 24} Let $\fe$ be arbitrary (not necessarily nondegenerate) continuous deformation of a conformal map $f(t)=t+2\pi it^{k+1}(1+O(t))$, 
$\a_i(\var)$, $i=0,\dots,k$, be its continuous fixed point families (not necessarily distinct), $\a(\var)=\a_j(\var)$ be one of them. There  exists a unique continuous in $\var$ (including $\var=0$) polynomial vector field family  $$w_\var(t)=a(\var)\pr$$ 
of degree $k+1$ with singularities at the fixed points of $\fe$ such that 
$$\fe\circ g^{-1}_{w_\var}(t)=t+O_1(t,\var),\tag5$$ 
$$O_1(t,\var)=O(\we(t)(t-\a(\var))), \ \ 
\frac{dO_1}{dt}(t,\var)=O(\we(t)(t-\a(\var)))+O(\frac{\partial(\we(t)(t-\a(\var)))}{\partial t}),$$ 
as $t,\var\to0$. (Recall (Subsection 2.4) that for a vector field $v$ by $g^s_v$ we denote its time $s$ flow map.)
 This follows from continuity of the family $\fe$ and is proved by straightforward calculation of the coefficient family  $a(\var)$ we are looking for: \  $a(\var)=
 {\prod_{s\neq j}(\a(\var)-\a_s(\var))}^{-1}\ln(1+2\pi i\prod_{s\neq j}(\a(\var)-\a_s(\var)))$, in particular, $a(0)=2\pi i$. In the case, when the family $\fe$ is nondegenerate, the 1-jet of the perturbed field $\we$ at its singular point $\a(\var)$ coincides with that of the canonic generator $\ve$ of the perturbed map $\fe$. 
\endremark

The previous inequality 2) holds a  priori in a neighborhood of $\a(\var)$ (that depends on $\var$) with $b=4\pi$ for small $\var$ 
(the last statement of the previous Remark).  To show that it  holds in a large domain, we construct a uniformly bounded family $G_\var$ of positive continuous functions   
in appropriate $\fe^{-1}$- invariant domains $D_\var\ni\a(\var)$, $G_\var(\a(\var))=2$,   such that the tangent disc field 
$T_G=(|\dot t|<G_\var(t)|\we(t)|)$ is $\fe$- invariant in $D_\var$. Then the 
inclusion $\ve\in T_G$, which is valid in a neighborhood of $\a(\var)$ by construction and the last statement of the previous Remark, will be also valid in the orbits of $\fe$ in $D_\var$ that go 
from $\a(\var)$ ($\fe$- invariance of $\ve$ and $T_G$).  
Therefore, this inclusion will hold in the whole domain $D_\var$: the $\fe^{-1}$- orbit of each point in $D_\var$ converges to $\a(\var)$, by $\fe^{-1}$- invariance of $D_\var$ and lemma 2.2 in \cite{10} (page 74). This together with the uniform boundedness of the functions $G_\var$ will prove the inequality 2) for appropriate $b$ independent on $\var$. 
 
For the construction of the domains $D_\var$ and the functions $G_\var$ we 
consider the vector fields $\we$ and a continuous in $\var$ (including $\var=0$) family $w_{0,\var}=e^{i\theta(\var)}\we$, $\theta(\var)\in\Bbb R$, of their constant multiples with the same moduli such that the perturbed fields $w_{0,\var}$ (correspondent to nonzero parameter values) have purely imaginary multipliers at their  
singularities $\a(\var)$ (so, the latters are centers). The possibility of choice of such a family continuous at $\var=0$ follows from the convergence of the arguments of the multipliers at $\a(\var)$ of the fields $\we$. The latter follows from the convergence assumption for $\arg\a(\var)$ in the case, when $k=1$, and the definition of nondegeneracy in the case, when $k\geq2$. We prove the statements from the item preceding Remark 24 for the following domains $D_\var$. 

\definition{Definition of the domains $D_\var$} Let $w_{0,\var}$ be a degree $k+1$ polynomial vector field family in $\Bbb C$ dependent continuously on the parameter $\var$ (including $\var=0$) with singularity 
families converging to 0, as $\var\to0$. Let $\a(\var)$ be its continuous singularity family with purely imaginary (nonzero) multiplier for all $\var\neq0$. Let $U$ be a disc centered at zero in the phase plane. For any $\var\neq0$ define $R_\var$ to be the rotation basin of 
the singular point $\a(\var)$, which is the domain saturated by the closed trajectories of the field $w_{0,\var}$ surrounding $\a(\var)$. Define $D_\var$ to be the maximal subdomain in $R_\var\cap U$ bounded by closed trajectory.
\enddefinition

\example{Example 4} In the conditions of the previous Definition let $k=1$ (i.e., the vector field family $w_{0,\var}$ be quadratic). Then all the closed trajectories of the perturbed field $w_{0,\var}$ (and in particular the boundary of $D_\var$) are circles symmetric with respect to the line passing through its  singularities.   
The rotation basin $R_\var$ of the singular point $\a(\var)$ is the half-plane $V_\var$ containing $\a(\var)$ and  bounded 
by the line orthogonal and bisecting the segment joining the singularities. Thus, $D_\var\subset V_\var$ for any $\var\neq0$. The line bounding $V_\var$ has limit, as $\var\to0$, which follows from the continuity of the family $w_{0,\var}$ at $\var=0$ and the center assumption on $\a(\var)$. The disc  $D_\var$ tends to the disc $D_0$ contained in the limit half-plane 
$V=\Lim_{\var\to0}V_\var$. Its boundary $\partial D_0$  is tangent to the boundary of $V$  at 0 and to that of $U$ (see Fig.5a,b). In particular, $D_\var$ satisfy the statement 3) following  the Definition of the domains $\ov$ from the beginning of the Subsection. 
\endexample

\hskip-2cm  {\bf Fig.5a,b}

Let us construct the function $G_\var$. Let $R_\var$ be the rotation basin from the previous Definition. For $t\in R_\var$ define $L(t)$ to be the (closed) trajectory of the field $w_{0,\var}$ passing through $t$. Define the function 
$G'_\var(t)$ in $t\in R_\var$ to be equal to the maximal distance between $\a(\var)$ and a point of $L(t)$. Let $N>0$. Put $G_\var=2+NG'_\var$. We show that for $N$ large enough the correspondent function family $G_\var$ is a one we are looking for (i.e., the correspondent tangent disc field $T_G:\ (|\dot t|<G_\var|\we|(t))$ is $\fe$- invariant in $D_\var$, whenever $U$ and $\var$ are small enough).  In the proof 
of this statement we use the following properties of the fields $\we$ and $w_{0,\var}$. 

\remark{Remark 25} 
Let $\fe$ be a nondegenerate family, $\a(\var)$ be its repelling fixed point family, 
$\we$, $w_{0,\var}$ be the correspondent vector field families from the item 
preceding the previous Definition. The fields $\we$ and $w_{0,\var}$ commute. The angle between 
them is constant in $t$ and bounded away from $\pi\Bbb Z$ uniformly in $\var$. 
The field $\we$ has repelling singular point at $\a(\var)$. It is directed outside the domains bounded by closed trajectories of the field 
$w_{0,\var}$ surrounding $\a(\var)$, so 
these domains (and in particular, $D_\var$) are $-\we$- invariant. These  statements follow from nondegeneracy of $\fe$ and the definitions of the families $\we$ and $w_{0,\var}$. The correspondent function $G_\var'$ (and hence, $G_\var$ also) from the last item increases along the trajectories of the perturbed field $\we$. 
In particular, the correspondent tangent disc field $T_G$ from the same item 
is $\we$- invariant (the last statement and $\we$- invariance of the field $\we$). 
\endremark 
\define\dist{\operatorname{dist}}

Firstly we prove Theorem 7 in the case, when $k=1$. For larger $k$ its proof is 
analogous and will be discussed at the end of the Subsection.

 Let us prove that the domains $D_\var$ are $\fe^{-1}$- invariant, 
and moreover, so is any subdomain in $D_\var$ bounded by closed trajectory of the field $w_{0,\var}$ (or equivalently, the 
function $G_\var'$ increases along the orbits of $\fe$ in $D_\var$) whenever $U$ and $\var$ are small enough. Moreover, we prove the following estimate of growth of the function $G_\var'$ along the $\fe$- orbits: there exists a $b>0$ such that the inequality 
$$G_\var'(t)-G_\var'(\fe^{-1}(t))\geq b|\we(t)|\tag6$$ 
holds in $D_\var$ whenever $U$ and $\var$ small enough. Further we use (6) in the proof of existence of invariant tangent disc field.    

Firstly let us prove (6) for the maps $g_{\we}^{-\frac12}$ instead of $\fe^{-1}$. The function $G_\var'$ decreases along the orbits of $g_{\we}^{-\frac12}$ (the previous Remark). We use the following estimate of the distance of a point $t$ to the $g_{\we}^{-\frac12}$- image 
of the correspondent curve $L(t)$: for any bounded neighborhood $U$ of zero in the $t$- line there exists a $b>0$ such that for any $\var$ small enough and $t\in R_\var\cap U$ 

$$\dist(t,g_{\we}^{-\frac12}(L(t)))\geq b|\we(t)|.\tag6'$$
Let us prove (6'). By $\beta(\var)$ denote the angle between $\we$ and $w_{0,\var}$ (which is bounded away from $\pi\Bbb Z$ by the previous Remark). 
We use the fact that the domain $M(t,\var)=\{ g^s_{\we}(t), \ \ |s|<\frac12 \sin\beta(\var)\}$ is disjoint from the curve $g^{-\frac12}_{\we}(L(t))$. This follows from the definition of $\beta(\var)$ and the fact that the fields $\we$ and $w_{0,\var}$ commute. (In fact, the boundary of the domain $M(t,\var)$ touches this curve.) Let $0<b<\frac12\inf\sin\beta(\var)$. The domain $M(t,\var)$ contains $b|\we(t)|$- neighborhood of $t$, whenever $\var$ and $t$ are small enough, since $g^s_{\we}(t)-t=s\we(t)(1+o(1))$, as $|s|\leq1$ and $t,\var\to0$. This together with the previous statement proves (6').
  
Now let us prove inequality (6) for the maps $g_{\we}^{-\frac12}$. Recall that the restriction of the function $G_\var'$ to the trajectories of the field $w_{0,\var}$ is constant. Let $t\in R_\var$. The curves $L(t)$ and  $g^{-\frac12}_{\we}(L(t))$ are circles surrounding $\a(\var)$ symmetric with respect to the line passing through the singularities of the field $w_{0,\var}$ (Example 4). Let $\tau$ be the point of the second curve having maximal distance to $\a(\var)$. Then $\tau$ lies in the same line and is separated from 0 by $\a(\var)$ in the line. By definition, $G_\var'|_{g_{\we}^{-\frac12}(L(t))}\equiv|\tau-\a(\var)|$.   Let $\tau'$ be the point of the curve $L(t)$ that lies in this line on the same side from $\a(\var)$, as $\tau$. Then $\tau'$ is the point of the curve $L(t)$ of maximal distance to $\a(\var)$ (so, $G_\var'(t)=|\tau'-\a(\var)|$). Therefore, 
$$G_\var'(t)-G_\var'(g_{\we}^{-\frac12}(t))=|\tau'-\a(\var)|- |\tau-\a(\var)|=|\tau'-\tau|\geq b|\we(\tau')|.$$ 
The last inequality holds whenever the curve $L(t)$ is close enough to 0 and $\var$ is small enough, i.e., the initial point $t$ lies in the domain $D_\var$ correspondent to $U$ and $\var$ small enough (estimate (6')). Now estimate (6) for the maps  $g_{\we}^{-\frac12}$ follows from the last inequality and the fact that the restriction $|\we||_{L(t)}$ takes its maximal value at $\tau'$ (since so does $|p(t,\var)||_{L(t)}$ by definition).

Now let us prove (6) for the maps $\fe^{-1}$. To do this, we use the estimate (6) for 
the maps $g^{-\frac12}_{\we}$ proved before, estimate (6') for the maps $g_{\we}^{\frac12}$ (which is valid, in fact, for any fixed real time $\we$- flow map and is proved in the same way, as for the time value $-\frac12$ above) and the asymptotic formula $\fe^{-1}(t)-g_{\we}^{-1}(t)=o(\we(t))=o(\we(g_{\we}^{-1}(t)))$, as $t,\var\to0$ (formula (5)). Then for any $U$ and $\var$ small enough for any $t\in D_\var$ the image $\fe^{-1}(L(t))$ is contained in the domain bounded by the curve 
$g_{\we}^{-\frac12}(L(t))$. This follows from  estimate (6') applied to the map $g_{\we}^{\frac12}$ at the point $g_{\we}^{-1}(t)$ and the last asymptotic formula. Then 
$$G_\var'(t)-G_\var'(\fe^{-1}(t))\geq G_\var'(t)-G_\var'(g_{\we}^{-\frac12}(t)).$$ 
This together with the inequality (6) for the maps $g_{\we}^{-\frac12}$ proves 
the same inequality for the maps $\fe$. Estimate (6) is proved.

Now let us show that for any $U$, $\var$ small enough and $N$ large enough the tangent disc field $T_G$ is $\fe$- invariant in $D_\var$.
To do this, we use (6) and the following estimate of the difference between the field $\we$ and its $\fe$- image: 
there exists a $d>0$ such that for any $U$, $\var$ small enough the inequality  
$$|(\fe)_*\we-\we|\leq d|\we|^2\tag7$$
holds in $D_\var$. In the proof of this estimate we use the following statement: 

 \item{$(*)$} for any $\var$ and $t\in R_\var$ \ $\a(\var)$ is the closest to $t$ of the singularities of the field $\we$.

By (5) and $g_{\we}^{-1}$- invariance of $\we$, 
$$\aligned
& ((\fe)_*\we-\we)(t)=((\fe\circ g_{\we}^{-1})_*\we-\we)(t)\\
& = O((\we(t))^2)+ O(\we(t)\frac{\partial((t-\a(\var))\we(t))}{\partial t}).
\endaligned$$ 
The last term in the right-hand side of this formula is $O((\we(t))^2)+ O((t-\a(\var))^2\we(t))$, as $t,\var\to0$. The second term in the last expression is $O((\we(t))^2)$, as $t\in R_\var$, $\var\to0$ (statement $(*)$). This proves (7). 

For the proof of invariance of the disc field $T_G$ for large $N$ let us 
calculate and estimate in $D_\var$ the difference of the radius of disc of the field $T_G$ and that of disc of its image 
$(\fe)_*T_G$. By definition, this difference is equal to 
$$\aligned
& (2+NG_\var')|\we|-(2+NG_\var'\circ\fe^{-1})|(\fe)_*\we|\\
& =N(G_\var'-G_\var'\circ\fe^{-1})|\we|-(2+NG_\var'\circ\fe^{-1})
(|(\fe)_*\we|-|\we|).\endaligned\tag8$$

It suffices to show that this difference is positive in $D_\var$, whenever $N$ is large enough and $U$, $\var$ are small enough. 
The first term of the difference in the right-hand side of (8) is 
greater than $Nb|\we|^2$ by (6). The module of the second term is not greater than $d|\we|^2(2+NG_\var'\circ\fe^{-1})$ by (7). The latter can be made arbitrarily smaller with respect to the former ($Nb|\we|^2$) (and hence, to the first term in the right-hand side of (8)) in $D_\var$, if $N$ be chosen large enough and $U$, $\var$ small enough. Indeed the first term in its brackets 
is a constant  independent on $N$ and $\var$; the second term is less than $N$ times the radius of the disc $U$ (which is greater than  $G_\var'|_{D_\var}$ by definition).  Thus, the second term 
in (8) can be made less than $Nb|\we|^2$, then the difference  (8) will 
be positive. The invariance of the tangent disc field is proved.

Now for the proof of Theorem 7 in the case, when $k=1$, it suffices to prove the  statement 1) following the Definition of the domains $\ov$ from the beginning of the Subsection. Let $V$ be the limit half-plane from Example 4. The statement 1) holds for any sector  $C$ whose closure 
punctured at 0 is contained in $V$. The sector $V$ is good: it contains the unique imaginary dividing ray, the latter lies in $S$ by definition. Therefore, the sector $C$ as above may be chosen to contain this ray as well and to be contained in $S$ (then it will be good). Theorem 7 is proved in the case, when $k=1$. 

Now let us prove Theorem 7 in the case, when $k\geq2$. Its proof in the previous case remains valid for larger $k$ modulo the statements 1) and 3) following the Definition of the domains $\ov$ from the beginning of the Subsection,  estimate (6) and statement $(*)$ from the item following (7). In the proof of these statements we use the following

\proclaim{Proposition 3} Let $k\geq2$, $w_{0,\var}$ be a family of degree $k+1$ 
polynomial vector fields in complex line with the coordinate $t$ depending continuously on the parameter $\var$ such that the (nonperturbed) field $w_{0,0}$ is a constant multiple of the monomial $t^{k+1}$ and the perturbed field (correspondent to nonzero parameter value) 
has $k+1$ singularities that form continuous families $\a_i(\var)$, $i=0,\dots,k$, satisfying the  asymptotic polygon regularity statement from Remark 10. Let $\Delta$ be the correspondent asymptotic regular polygon. Let $\a(\var)$ be a center singularity family (i.e., the correspondent multipliers be purely imaginary whenever $\var\neq0$), $A$ be the correspondent vertex of $\Delta$. Let $V$ ($V_\var$) be the radial sector with the angle $\frac{\pi}k$ containing $A$ and symmetric with respect to the radial line of $A$ (respectively, $\a(\var)$). Let $W_\var$ be the sector 
containing $\a(\var)$ bounded by the radial rays orthogonal to the singularity polygon sides neighboring at $\a(\var)$. For any $\var$ small enough the rotation basin $R_\var$ of the center singularity $\a(\var)$ is contained in $W_\var$ and bounded by the  trajectory of the field $w_{0,\var}$ that goes from infinity and returns to it  in the asymptotic directions of the boundary rays of $V_\var$ (not necessarily approaching the latters). The rotation basin family $R_\var$ converges to $V$, 
as $\var\to0$ (see Fig.6 a,b in the case, when $k=2$). For any neighborhood $U$ of zero in the phase line there exists an $r'>0$ such that the correspondent domains $D_\var$ from the previous Definition satisfy the statement 3) following the Definition of the domains $\ov$ at the beginning of the Subsection. 
\endproclaim

\hskip-2cm Fig. 6a,b

\proclaim{Corollary 5} In the conditions of the previous Proposition let $U$ be a neighborhood of zero in the $t$- line, $D_\var$ be the correspondent family 
of domains from the previous Definition. The family $D_\var$ converges to a subdomain $D_0\subset U$ bounded by the trajectory of the field $w_{0,0}$  that goes from zero and returns to it with asymptotic tangency to the boundary rays of $V$ and touches the boundary of $U$ (see Fig.6b in the case, when $k=2$).
\endproclaim 

\demo{Proof of Proposition 3} The proofs of the statements of Proposition 3 on the boundary of  $R_\var$ and the convergence of the latter are implicitly contained in subsection 6 of \cite{6}. Let us show that for any $\var$ small enough $R_\var\subset W_\var$. To do this, we show that 
for any $\var$ small enough the field $w_{0,\var}$ is transversal to the boundary rays of $W_\var$. Indeed in the case, when the singularity polygon is regular and centered at 0, the field $w_{0,\var}$ has constant nonzero angles  with the boundary rays of $W_\var$: the correspondent real tangent line field is constant in the union of the symmetry lines of the polygon (subsection 6 in \cite{6}); therefore, it is orthogonal to the radial ray of the singular point $\a(\var)$ (the latter is center); hence, its angles with the boundary rays of the sector $W_\var$ are constant and equal to $\frac{\pi}2-\frac{\pi}{k+1}>0$. Now for the proof of transversality in the general 
case we consider the new vector field family $\tilde\we$ obtained from $\we$ by the radial homothethy family from Remark 10, which transforms the singularity polygon of the latter to a polygon with unit diameter converging to the correspondent  regular polygon $\Delta$.  
The transversality statement from the beginning of the item is equivalent to that for the new family $\tilde\we$. This statement for $\tilde\we$ follows from continuity in $\var$ of the correspondent real tangent line field  
and the same statement for the limit field $\tilde w_0$, which has a regular  singularity polygon. 

We prove the inclusion $R_\var\subset W_\var$ by contradiction, using only the previous transversality statement. Suppose $R_\var\not\subset W_\var$. Then there exists a closed trajectory of the field 
$w_{0,\var}$ in $R_\var$ tangent to a boundary ray of $W_\var$, since, by assumption, 
there exists at least one closed trajectory in $W_\var$ ($\a(\var)\in W_\var$). 
This contradicts the transversality statement. The inclusion statement of Proposition 3 is proved. 

Now let us prove the last statement of Proposition 3. It is reduced to its particular case, when the singularity polygons are regular, analogously to the discussion from the item following Proposition 3. In this case for all $\var$ the domain $D_\var$ is symmetric with respect to the radial line of $\a(\var)$. Its intersection with this line is connected: otherwise $D_\var$, which is simply connected by definition, would not be so 
by symmetry. The last statement of Proposition 3 follows from the convergence of this intersection to that of the limit line with the limit domain $D_0$ from Corollary 5, which is a nonempty radial interval of the limit of the radial ray of $\a(\var)$. Proposition 3 is proved.   
\enddemo

Now let us prove the four statements from the item preceding Proposition 3. The statement 3) follows from the last statement of Proposition 3. Let us prove the statement 1). The inclusion $D_0\supset C^r$ holds for any sector $C$ whose closure with zero deleted is contained in the sector $V$ from Proposition 3 and appropriate $r>0$ depending on $C$ (Corollary 5). By definition (Remark 11), the sector $V$ is contained in $S$ and is good (contains the same imaginary dividing ray, as $S$). Therefore, the sector $C$ as above may be chosen to be good as well. 

Let us prove statement $(*)$. Let us consider that $\a(\var)$ is the singularity closest to 0 for all $\var$: one can achieve this by applying appropriate continuous translation family in the $t$- line so that the asymptotic regularity statement from Remark 10 will remain valid. Then $(*)$ follows from Proposition 3 (the inclusion $R_\var\subset W_\var$ and the definition of the sector $W_\var$).

Let us prove estimate (6). Its prove is done in the same way, as in the previous case modulo the same inequality for the maps $g^{-\frac12}_{\we}$ instead of $\fe^{-1}$. Let us prove this inequality. To do this, we use estimate (6'), which is valid for any fixed real time flow map of the field $\we$ and is proof repeats that in the previous case. For $t\in R_\var$ by $\tau(t)$ ($\tau'(t)$)  denote the point of the curve $g^{-\frac12}_{\we}(L(t))$ (respectively, $L(t)$)  that has the maximal distance to $\a(\var)$. By definition, the functions $\tau'(t)$ and $\tau(t)$ are constant in the trajectories of the field $w_{0,\var}$. Then 
$$G_\var'|_{L(t)}-G_\var'|_{g^{-\frac12}_{\we}(L(t))}\geq\dist(\tau(t),\ L(t))\geq b|\we(\tau(t))|$$
(the first inequality follows from the definition of the functions  $G_\var'$,  $\tau$ and $\tau'$; the last inequality follows from the estimate (6') applied to the map $g^{\frac12}_{\we}$ at the point $\tau(t)$). Now for the proof of (6) for the maps $g^{-\frac12}_{\we}$ it suffices to show that there exists a $c>0$ such that for any $U$, $\var$ small enough and $t\in D_\var$ 
$$|\we(t)|\leq c|\we(\tau(t))|.\tag9$$  
In the proof of the last inequality we use statement $(*)$ proved before and the fact that for any $U$, $\var$ small enough and $t\in D_\var$  $$|t-\a(\var)|\leq2|\tau(t)-\a(\var)|=2G_\var'|_{g^{-\frac12}_{\we}(L(t))}
.\tag10$$
 The latter follows from  the inequality $|t-\a(\var)|\leq G_\var'|_{L(t)}$ and the asymptotic formula $G_\var'|_{L(t)}=G_\var'|_{g^{-\frac12}_{\we}(L(t))}(1+o(1))$, as $\var, L(t)\to0$ (the asymptotic formula $t'=g^{-\frac12}_{\we}(t')(1+o(1))$, as $\var,t'\to0$). 

For the proof of (9) it suffices to show that for any other singularity family  $\a_i\neq\a$ $|t-\a_i(\var)|\leq4|\tau(t)-\a_i(\var)|$, whenever $U$, $\var$ are small enough and $t\in D_\var$. This is implied by the following inequalities:
$$|t-\a_i|\leq |t-\a|+|\a-\a_i|; \ |t-\a|\leq2|\tau(t)-\a|<2|\tau(t)-\a_i|;$$ 
$$|\a-\a_i|\leq|\tau(t)-\a|+|\tau(t)-\a_i|<2|\tau(t)-\a_i|.$$
The inequalities in the first line follow from the triangle inequality, (10) and statement (*)  
respectively. Those in the second line follow from the triangle inequality and (*) respectively. Estimate (9) (and hence, (6) also) is proved. The proof of Theorem 7 is completed. 

\subhead 3.5. Proof of the two last statements of Proposition 2 \endsubhead

In the conditions of Proposition 2 consider the new family $w_{0,\var}=e^{i\theta(\var)}w(\var)$, $\theta(\var)\in\Bbb R$, of multiples of the fields $w(\var)$ having purely imaginary multipliers at $\a(\var)$. Without loss of generality we consider that the argument of the singularity $\a(\var)$ has limit, as $\var\to0$. In the case, when $k\geq2$, this condition holds by definition. In the case, when $k=2$, one can achieve this by taking appropriate subsequence of arbitrary convergent sequence $\var_n\to0$ instead of all the parameter space. Then the family $w_{0,\var}$ can be chosen to be continuous at $\var=0$. Let $D_\var$ be the family of domains from the previous Definition correspondent to the neighborhood $U_t$. Then $D_\var\subset\ov$ by definition. This together with the last statement of Proposition 3 proves the last statement of Proposition 2. Let $k\geq2$. Let us prove that the connected component of the intersection $\ov\cap\Omega$ containing $\a(\var)$ tends to $\Omega$, as $\var\to0$. This will prove Proposition 2. To do this, it suffices to show that just some subdomain of the intersection converges to $\Omega$: the last statement of Proposition 2 will imply that this holds for the intersection component containing $\a(\var)$, as in Subsection 3.2. Let $D_0$ be the limit domain from the previous Corollary. Then $D_0\subset\Omega$. Let $C$ be a good sector, $r>0$, such that $C^r\subset D_0$. We consider that $C$ contains the repelling ray correspondent to the limit vector field $w(0)$ (one can achieve this by choosing $C$ large enough). Then a connected component of the intersection $D_\var\cap C^r$ (and hence, $\ov\cap C^r$) converges to $C^r$. Consider the family of domains saturated by the trajectories of the fields $w(\var)$ in $\Omega$  that start in  $C^r\cap\ov$. This family converges to $\Omega$ by definition.  The proof of the two last statements of Proposition 2 is completed.  

\subhead 3.6. Another proof of Theorem 3(4) that uses Theorem 7\endsubhead 

Let us sketch another proof of Theorems 3 and 4. Let $D=\{|t|<\delta, z=\delta\}$ 
be  a disc transversal to the phase curves of the vector fields from the 
family $\le$. Consider the monodromy map family $\fe$ of the fields $\le$ in $D$  correspondent to going around the singular point $0$ in the separatrix $t=0$ 
of the nonperturbed field. 
The family $\fe$ is nondegenerate in the sense of the Definition from 
Subsection 2.4.C. (In fact, the nondegeneracy of the monodromy family 
is equivalent to that of the family $\le$ under considerations.) 

As it is shown below, Theorems 3 and 4 are implied by Theorem 7 
and the two following  Remarks. 

\remark{Remark 26} Let (1) be a vector field as in Section 1 holomorphic 
in a neighborhood $U=\{|z|<2\delta\}\times\{|t|<\delta\}$ of zero.  
 Let $D=\{ z=\delta\}\cap U$ be a disc transversal to (1) in the line $z=\delta$   equipped 
with the coordinate $t$. Then the monodromy map of (1) in $D$ correspondent to  the 
counterclockwise going around zero in the separatrix $t=0$ is holomorphic in 
a neighborhood of the point $t=0$ in $D$ and satisfies the conditions of  
 Theorem 5. Let $S$ be a good sector in the $t$- axis, $v$ be the correspondent 
sectorial canonic generator of the monodromy map, $\tau$ be its
 complex time. The correspondent 
sectorial canonic first integral of (1) restricted to $D$ is a constant 
multiple of the function $e^{2\pi i\tau}$. This follows from the same statement  for the formal normal form $(1)_n$ and Theorem 1.
\endremark

\remark{Remark 27} Let 
$$\cases
& \dot z=\nu z\\ 
& \dot t=\mu t
\endcases$$
be a linear vector field with $\frac{\mu}{\nu}\notin\Bbb R$,  $|\frac{\mu}{\nu}|<\frac12$. Let $D=\{|t|<\delta,\ z=\delta\}$ be a disc transversal to the field equipped with the coordinate $t$. Let $M$ be 
 the monodromy map of the field in $D$ correspondent to going around zero in 
the separatrix $t=0$. The multiplier of $M$ at its fixed point $t=0$ is equal to  $\la=e^{2\pi i\frac{\mu}{\nu}}$. The 
canonic first integral of the vector field and the 
canonic generator of $M$ at the fixed point are well-defined   
($|\la|\neq1$). Let $\tau$ be the (multivalued) complex time in $D$ correspondent to 
the generator. The restriction of the canonic integral to $D$ (which is a multivalued function branched at $t=0$) is a constant multiple of 
the function $e^{2\pi i \tau}$.
 This statement remains valid for a nonlinear vector field having a singular point with multipliers as at the beginning of the item. This follows from Theorem 2. 
\endremark

Now let us sketch the proof of the convergence statement of Theorem 3(4). 
The restrictions of  appropriately normalized canonic first integrals 
of $\le$ to the transversal disc $D$ 
converge to that of the correspondent sectorial canonic first integral. This 
follows from the two previous Remarks and  Theorem 7 applied to the monodromy 
family $\fe$. Now for the proof of Theorem 3(4) it suffices to show that  the 
canonic integrals of the perturbed fields extend analytically to domains satisfying statement 1) of Theorem 3 (along the phase curves, where they are 
constant), and this extension depends continuously on the parameter.
To do this, one shows that for appropriately chosen neighborhoods 
$U_z$, $U_t$ of zero in the coordinate lines and appropriate domains $\ov\subset U_t$ (e.g.,  the subdomains $\ov^2$ from Subsection 3.3) each point of the  domain $\widetilde{\ov}=U_z\times\ov$ can be connected in the latter to a point of $D$  where the 
monodromy generator is well-defined by path lying in the phase curve of $\le$, 
and one can choose these paths to be dependent continuously on the parameter. 
The details of the proof are omitted to save the space. 

 \subhead 3.7. Proof of Lemma 1\endsubhead
 \define\dve{(2)_{\var}}
  \define\tre{(3)_{\var}}
 \define\lje{\lambda_i(\var)}
\define\mje{\mu_i(\var)}
\define\Lje{L_i(\var)}
\define\fje{F_{i,\var}}
\redefine\aje{\a_i(\var)}
 
Let 
 $$\hskip-4.8cm\dve\hskip4.8cm\cases 
 & \dot z=f(z,t,\var)\\
 & \dot t=g(z,t,\var)\endcases$$
 be arbitrary continuous deformation of a vector field (1). This means that 
 the functions $f$ and $g$ are continuous in $\var$ and holomorphic in $(z,t)$, 
$f(z,t,0)=z+O(|z|^2+|t|^{k+1})$, $g(z,t,0)=t^{k+1}$. We show that $\dve$ 
is orbitally analytically equivalent to a family $\le$, i.e., phase curves 
of the fields $\dve$ can be transformed to those of $\le$ by continuous family 
of analytic changes of the variables $(z,t)$ defined for all the parameter 
values small enough. This will prove Lemma 1.

Without loss of generality we assume that all the singular points of the 
vector fields $\dve$ lie in the $t$- axis. One can achieve this by applying 
the family of coordinate changes of the variable $z$ defined by the formula 
$\tilde z=f(z,t,\var)$. Then the family of vector fields under consideration is 
of the type

 $$\hskip-2cm\tre\hskip2cm\cases
 & \dot z=z(1+q(z,t,\var))+g(t,\var)\pr\\
 & \dot t=\pr+zR(z,t,\var)\endcases$$ (up to multiplication by family of 
 nonzero holomorphic functions), 
where the functions $q$, $g$, $R$ are holomorphic in $(z,t)$ and continuous in the parameter $\var$, as are the constants $\a_i(\var)$, $\a_i(0)=0$, 
$q(0,0,0)=0$, $R(z,t,0)\equiv0$. The singular points of these vector fields are $(0,\aje)$. 
Let us show that there exists a family of changes of the variable $t$ that transforms 
$\tre$ to an analogous family with $R\equiv0$. This will prove Lemma 1.  

Let $\lje$ be the family of eigenvalues of the linearization operator of 
$\tre$ at the continuous singular point family $(0,\aje)$ such that the  correspondent eigenline family approaches the line tangent to the $z$- axis, as $\var\to0$ (so, 
$\lambda_i(0)=1$). Let $\mje$ be the other eigenvalue family. 
Then $\mu_i(0)=0$. For any $\var$ small enough 
there exists a unique separatrix $\Lje$ passing through the singular point $(0,\aje)$  
and tangent to the eigenline of the linearization 
operator at $(0,\aje)$ correspondent to the eigenvalue $\lje$. In general, for any planar vector field in a neighborhood of its singular point with multipliers $\lambda$ and $\mu$, $|\lambda|>|\mu|$, there exists a unique regular holomorphic curve tangent to the field and passing through the singularity with tangency to the $\lambda$- eigenline of the linearization operator. In the case, when $\frac{\mu}{\lambda}\notin\Bbb R_+$, this follows from either Theorem 2 (when this ratio is not real), or the analytic version of Hadamard - Perron theorem \cite{12} (when it is either real negative, or zero). In the case, when this ratio is positive,  the existence and uniqueness of separatrix statement is reduced to the previous one by applying blowing-up. The last 
statement from the last item is equivalent to the existence of a $t$- coordinate  change family that rectifies all these separatrices simultaneously. To prove 
this its reformulation, we use the convergence of the separatrices of the perturbed field to that of the nonperturbed field. This is implied by the following  

\proclaim{Proposition 4} Let $\tre$ be as at the beginning of the Subsection and in the item following its formula, $\Lje$ be a separatrix family from the last item. There exist a neighborhood $U=U_z\times U_t$ of zero in the phase space such that for 
any $\var$ small enough the local separatrices $\Lje\cap U$ 
are graphs $t=\fje(z)$ of functions $\fje(z)$ holomorphic in $U_z$
and depending continuously  on the parameter (including $\var=0$, where $F_{i,0}\equiv0$) such that $|(\fje)'_z|\leq1$, $|\fje(z)|\leq|z|$. 
\endproclaim

\redefine\lje{\Lje}
\define\loe{L_0(\var)}
\define\prl{\prod_{i=1}^k(t-\aje)}
\define\cre{(4)_{\var}}
\define\fle{F_{1,\var}}
\define\ale{\a_1(\var)}
\define\lle{L_1(\var)}

\demo{Proof} Let $(0,\aje)$ be a continuous family of singularities of $\tre$.  
Let $K$ and $\widetilde K$ be respectively the tangent cone field 
and the cone  from the proof of Lemma 4 (Subsection 3.2). 
The inequalities of Proposition 4 are equivalent to the statement that 
the tangent lines to $\lje$ lie outside $K$ and $\lje$ lies outside 
$\widetilde K$. These inclusions hold a priori in a neighborhood of 
$(0,\a_i(\var))$ (that depends on $\var$). As in the proof of Lemma 2, 
to show that they hold in a large domain, 
we use the fact that there exists a neighborhood $U$ of zero where the field of  the complements to the cones $K$ (and hence, the complement to $\widetilde K$) are invariant with 
respect to the field (3)$_\var$ for all $\var$ small enough (proposition 2 in 
\cite{6}). Let us consider 
that  $U=U_z\times U_t$, where $U_z=\{|z|<\delta\}$, $U_t=\{|t|<2\delta\}$. 
Let $L_i'(\var)\subset\lje\cap U$ be the subdomain saturated by the trajectories of 
the field $\tre$ in $\lje\cap U$ that go from $(0,\a_i(\var))$. Then $TL_i'\cap K=\emptyset$, $L_i'\cap\widetilde K=\emptyset$, by invariance of the cone $K$ complement field.
Now for the proof of Proposition 4 (modulo continuity of the  
dependence on the parameter) it suffices to show that $L_i'(\var)$ is the graph  $t=\fje(z)$ of a holomorphic function  in $U_z$ (then $L_i'(\var)=L_i(\var)\cap U$). This is equivalent to the statement that 
$L_i'(\var)$ is 1-to-1 projected onto the disc $U_z$ in the $z$- line. This is implied  by the following properties of this projection: 1) it 
does not have critical points 
(the tangent lines to $L_i'(\var)$ are not parallel to the $t$- line, 
since they lie outside $K$);
2) it is a map "onto" $U_z$. The last statement follows from the fact that 
the trajectories that form $L_i'(\var)$ meet 
the boundary of $U$ at points whose projection images in the $z$- line lie in 
the boundary of the disc $U_z$. This follows from the inclusion $L_i'(\var)\subset U\setminus\widetilde K$ and the definition of $U$ (in the same way, as in the proof of lemma 1 from 
\cite{6}). The inequalities of Proposition 4 are proved. 
Let us prove the continuity of the family of the correspondent 
functions $\fje$ in $\var$, e.g., at $\var=0$ (for $\var\neq0$ the proof is analogous). The family $\fje$ converges to 0 uniformly in compact subsets of $U_z$, as $\var\to0$. This follows from its normality (inequalities of Proposition 4) and the uniqueness of the separatrix of the nonperturbed field transversal to the $t$- axis, as in the proof of Lemma 2 in Subsection 3.2. Proposition 4 is proved.
\enddemo

By Proposition 4, each family $\Lje$ is continuous in the parameter, and hence,  
can be rectified itself (i.e., transformed to a family of subdomains of complex lines parallel to 
the $z$- axis)  by a continuous family of changes of the variable $t$. Now 
we show that these families can be rectified simultaneously by $t$- variable changes depending continuously on the parameter. This will prove Lemma 1.
 
Without loss of generality we consider that the singular point family $(0,\a_0(\var))$ is identically zero, and the correspondent separatrices $\loe$  already lie in the $z$- axis. This means 
that the family $\tre$ under considerations is of the type 

 $$\hskip-2cm\cre\hskip2cm\cases 
 & \dot z=z(1+q(z,t,\var))+tg(t,\var)\prl\\
 & \dot t=t(\prl+zR(z,t,\var)),
 \endcases$$

\define\prd{\prod_{i=2}^k(t-\aje)}

Let us firstly show that there exists a continuous family of changes of the 
variable $t$ that preserves 
the line $t=0$ (in particular, the separatrix $\loe$) and the singularity family  $(0,\a_1(\var))$ and transforms the correspondent separatrix 
$\lle$ to a subdomain of the line $t=\a_1(\var)$. 
Let $\fle$ be the correspondent function from Proposition 4. 
The family of changes we are looking for is 
$\tau=t\frac{\ale}{\fle(z)}$. By construction, it preserves $\loe$ and 
rectifies $\lle$. Let us prove that these changes form a continuous family. 
To do this, it suffices to show that they converge to the identity, as 
$\var\to0$.   

For the proof of the last statement we show that 
$$(\ln\fle)'_z=\frac{(\fle)'_z}{\fle}\to0, \ \text{as} \ \var\to0$$
uniformly in compact sets. This together with the equality 
$\fle(0)=\ale$ will imply that $\frac{\fle(z)}{\ale}\to1$, and hence, the 
convergence of the coordinate changes under considerations to the identity.
The calculation of the derivative of the function $\fle$ yields 

$$\fle(z)'_z=\frac{dt}{dz}|_{\lle}=\frac{t((t-\ale)\prd+R(z,t,\var)z)}
{z(1+q(z,t,\var))+g(t,\var)t(t-\ale)\prod_{i=2}^k(t-\aje)}|_{t=\fle(z)}.$$
Therefore, 
$$(\ln\fle(z))'_z=\frac{\frac{(t-\ale)}z\prd+R(z,t,\var)}{1+q(z,t,\var)+g(t,\var)
\frac{t-\ale}zt\prd}|_{t=\fle(z)}.$$

The ratio $\frac{t-\ale}{z}|_{t=\fle(z)}$ in the right-hand side of the last formula has module at most 1 (the last inequality of Proposition 4). Therefore, 
there exists a neighborhood $U_z$ of zero in the $z$ -axis  (independent on $\var$) where the right-hand side has 
module less or equal to 
$2(\prod_{i=2}^k(|\fle(z)|+|\aje|)+|R(z,\fle(z),\var)|)$ for all $\var$  small enough. The last expression converges  
to 0, as $\var\to0$, since so do $\fle(z)$, $\aje$ and $R(z,t,\var)$. 
The statements on the convergence of the coordinate changes to 
the identity is proved. (This proves the simultaneous rectification statements 
for the separatrices in the case, when $k=1$.) 

We have already proved the possibility to rectify two separatrix families 
simultaneously by continuous family of changes of the $t$- coordinate. 
The possibility to rectify 3 or more separatrix families simultaneously  is  proved 
analogously by induction in their number. 
For example, let us show that if 2 separatrix families (say, $\loe$ and $\lle$) 
are already rectified (i.e., lie in lines $t=\a_0(\var)(=0)$, $t=\a_1(\var)$ respectively), then  
there exists a continuous family of changes of the 
$t$- coordinate that preserves them and rectifies another separatrix family (say, $L_2(\var)$). 

The rectifying family we are looking for is the family of changes 
$$\tau=t-t(t-\ale)\frac{F_{2,\var}(z)-\a_2(\var)}{F_{2,\var}(z)(F_{2,\var}(z)-\ale)}.$$
As before, it suffices to prove the convergence of this family to identity, as 
$\var\to0$. To do this in its turn, it suffices to show that 
$\ln(F_{2,\var})'_z=o(\a_2(\var)-\a_1(\var))$, as $\var\to0$. The proof of 
this statement repeats that of the analogous statement for $\fle$ proved 
before with obvious changes. Lemma 1 is proved.

\head 4. A higher-dimensional central manifold analogue of Theorems 3 and 4
\endhead

\subhead 4.1. Statement of results. Existence and uniqueness of sectorial central manifolds  of general higher-dimensional saddle-node vector fields. Their expression as limits  of separatrices of generic perturbation\endsubhead
\define\cn{\Bbb C^n}
\redefine\Re{\operatorname{Re}}

In this Subsection we consider saddle-node singularities of higher dimensional holomorphic vector fields, which are defined in exactly the same way, as 
in the two-dimensional case (at the beginning of the paper).

\remark{Remark 28} Any vector field with a saddle-node singularity is 
locally orbitally analytically equivalent to a field 
 of the type 
$$\cases
& \dot y=By+O(|y|^2+|t|^{k+1}) \\
& \dot t=t^{k+1}+O(|y|^2)\endcases\tag11$$
in a neighborhood of zero in $\cn$ with the coordinates $(y=(y_1,\dots,y_{n-1}), \ t)$, where $B$ is a nondegenerate upper-triangular block-diagonal matrix (each block corresponds to one eigenvalue). 
The reduction of arbitrary saddle-node vector field to this form is made 
by linear change of variables (to make the linearization matrix to consist of 
one-dimensional zero block (correspondent to the $t$- axis) and a matrix $B$ as above), subsequent killing the (nonresonant) monomials $t^ly_i\frac{\partial}{\partial t}$, $l<k+1$, and multiplication of the new field by appropriate nonzero holomorphic function.\endremark

Generically, a vector field (11) does not have a central manifold, i.e., a 
regular separatrix tangent to the $t$- axis at 0. At the same time, it always has 
a formal central manifold: there exists a unique formal $n-1$- component vector 
Taylor series 
$\hat q(t)$ without terms of degrees less than 2 such that the formal 
change $\tilde y=y-\hat q(t)$ of the variable $y$ transforms (11) to a (formal)  vector field tangent to the $t$- axis (i.e., having Taylor series that does not contain the terms $t^l\frac{\partial}{\partial y_j}$, which are nonresonant). 
Generically, this series diverges. At the same time there exists a 
covering $\bigcup_{j=0}^mS_j$ of a punctured neighborhood of zero in the $t$- 
plane by radial sectors such that in each sector $S_j$ there exists a unique $(n-1)$- dimensional holomorphic vector function 
$q_j(t)$ that is $C^{\infty}(\overline S_j)$, $q_j(0)=0$, whose graph $y=q_j(t)$ is tangent to the field. All the vector functions $q_j$ have the same asymptotic Taylor series  at zero coinciding with $\hat q$ (see Theorem 8 below).  

A radial sector $S$ satisfies all the statements from the last item, if 
it satisfies the conditions of the following 

\definition{Definition} Let (11) be a vector field as in Remark 28, $b_i$, $i=1,\dots,n-1$, be the eigenvalues of the matrix $B$. In this Subsection 
by {\it imaginary (real) dividing ray or line} correspondent to an eigenvalue $b_i$ we mean a radial ray or line from the set $\{ \Re(\frac{b_i}{t^k})=0\}$ 
(respectively, $\{\operatorname{Im}(\frac{b_i}{t^k})=0\}$). A radial sector 
in the $t$- plane is said to be good, if for any $b_i$ it contains exactly one correspondent imaginary dividing ray and its closure does not contain additional imaginary dividing rays.
\enddefinition

\example{Example 5} In the conditions of the last Definition let $n=2$,  $b_1\in\Bbb R$. In this case the new definitions of imaginary (real) dividing rays and good sectors  coincide with those from Section 2. 
\endexample

\proclaim{Theorem 8} Let (11) be a vector field as in Remark 28, $S$ be a good sector. There exists a unique $(n-1)$- dimensional vector function $q(t)$ holomorphic in a neighborhood of zero in $S$, continuous at 0  and having bounded derivative, $q(0)=0$, such that its graph $y=q(t)$ over $S$ is contained in a phase curve of (11). The vector function $q$ is $C^{\infty}(0)$ and has 
asymptotic Taylor series at 0. This series does not depend on the choice of the sector and coincides with the formal central manifold series $\hat q$ from the item following Remark 28. 
\endproclaim

In the case, when $n=2$, Theorem 8 is reduced to Theorem 1 (Subsection 2.1). Some higher-dimensional particular cases of Theorem 8 were proved in \cite{11}. Its proof in the general case is presented in Subsections 4.2-4.4. 
 
We consider a generic deformation 
$$\hskip-4cm(11)_\var\hskip4cm\cases
& \dot y=G(y,t,\var) \\
& \dot t=F(y,t,\var)\endcases$$ 
of a vector field (11) (correspondent to zero value of the deformation parameter $\var$) that in particular splits the degenerate 
singularity 0 of the nonperturbed field into $k+1$ distinct singularities of 
the perturbed field. The linearization operator of the perturbed field 
at each singularity has eigenvalue family that tends  to zero, as $\var\to0$. 
The singularity of a generically perturbed vector field possesses a unique local separatrix tangent to the eigenline correspondent to this eigenvalue. We show that for a generic deformation $(11)_\var$ these ("horizontal") separatrices of the perturbed  vector field converge to appropriate sectorial 
separatrices $y=q_j(t)$ of the nonperturbed field from Theorem 8. 

Before the statement of this result we firstly recall the Theorem on the existence of separatrix at a nondegenerate singularity of a holomorphic vector field. 

\proclaim{Theorem 9} Let a holomorphic vector field in $\cn$ have nondegenerate singularity at 0, $b_1,\dots,b_n$ be the eigenvalues of the correspondent linearization operator. Let $\frac{b_n}{b_i}\notin\Bbb R$ for any $i\neq n$. Let  $l_n$ be the eigenline correspondent to the eigenvalue $b_n$. Then 
there exists a unique separatrix of the field tangent to $l_n$. 
\endproclaim

Theorem 9 follows from theorem 4.3.2 in \cite{12} applied to a constant multiple of the vector field under considerations whose linearization operator eigenvalue  correspondent to $l_n$ is purely imaginary.
 
Without loss of generality we consider families  $(11)_\var$ such that all the singularities of the perturbed field lie in the $t$- axis (their $t$- coordinates will be referred to, as $\a_i(\var)$, $i=0,\dots,k$). One 
can achieve this by applying appropriate continuous family of changes 
of the variables $y$ (e.g., $\tilde y=G(y,t,\var)$).
Define $p(t,\var)=\prod_{i=0}^k(t-\a_i(\var))$. Then $G(0,\a_i(\var),\var)=0$, and hence, 
$G(y,t,\var)=By+o(|y|)+O(p(t,\var))$, as $(y,t,\var)\to0$. Without 
loss of generality we also consider that  $F(y,t,\var)=
p(t,\var)+O_2(y,t,\var)$, $O_2(y,t,\var)=o(|y|)$, as $(y,t,\var)\to0$  
(one can achieve this by applying appropriate continuous family of linear $t$- coordinate changes), and $\sum_{i=0}^k\a_i(\var)\equiv0$. 
 
We prove the statement from the item preceding Theorem 9 for the following  families  $(11)_\var$.

\definition{Definition} A family $(11)_\var$ as in the item following Theorem 9  is said to be nondegenerate, if each line passing through a singularity pair of the perturbed field intersects each real dividing line by angle bounded away from zero uniformly in $\var$, and in the case, when $k\geq2$, the correspondent polynomial $p(t,\var)$ satisfies the conditions from the second item of Subsection 2.3.B.
\enddefinition

 \remark{Remark 29} Let $(11)_\var$ be a nondegenerate family, $\a_j$ be its  continuous singularity coordinate family. In the case, when $k=1$, let $V_{j,\var}$ be the radial sector with the angle $\pi$ containing $\a_j(\var)$ and symmetric with 
respect to its radial line (see Fig.5a). The sectors $V_{j,\var}$ are good, 
and their boundaries do not accumulate to no imaginary dividing ray. The closure of their union is contained in a good sector. (The latter will be referred to, as $S_j$.) To the family $\a_j(\var)$ we put into  correspondence the sector $S_j$. In the case, when $k\geq2$, let $A_j$ be the vertex of the limit regular polygon from Remark 10 correspondent to $\a_j$, $V_j$ be the radial sector with the angle $\frac{\pi}k$ containing $A_j$ and symmetric with respect to its  radial line. The sector $V_j$ is good. To the family $\a_j(\var)$ we put into correspondence a good sector $S_j\supset V_j$. (The sectors $S_j$ do not cover 
a neighborhood of zero in the case, when $k\geq2$, cf. Remark 12.)
\endremark

\example{Example 6} In the case, when $n=2$ and the linearization operator 
of the nonperturbed field $(11)_0$ has real eigenvalues, the above nondegeneracy Definition and the definition of the sectors correspondent to the singularity coordinate families from the previous Remark coincide with those from Subsection 2.3 (with the only difference that in the case, when $k=1$, the sectors $S_j$ in Remark 29 are chosen to be large enough to contain the closure of the union (in $\var$) of the sectors $V_{j,\var}$). \endexample 

\remark{Remark 30} Let $(11)_\var$ be a nondegenerate family, $\a_j(\var)$ be its continuous singularity coordinate family. Let  
$\mu(\var)$ be the eigenvalue family of the correspondent linearization operator  that tends to 0, as $\var\to0$. For all $\var$ 
small enough the perturbed field $(11)_\var$ satisfies the conditions of 
Theorem 9 at $(0,\a_j(\var))$ with respect to the eigenvalue $b_n=\mu(\var)$. Moreover, for any other continuous eigenvalue family $b(\var)$ of the linearization operator  
the argument $\arg(\frac{\mu(\var)}{b(\var)})$ is bounded away from  $\pi\Bbb Z$ uniformly in all $\var$ small enough. In particular, for any small $\var$ 
the field $(11)_\var$ possesses a separatrix at $(0,\a_j(\var))$ (it will be referred to, as $\G_{j,\var}$) tangent to the eigenline 
of the linearization operator correspondent to $\mu(\var)$ (Theorem 9). (This eigenline approaches the line tangent to the $t$- axis, as $\var\to0$.) \endremark

\proclaim{Theorem 10} Let $(11)_\var$, $\a=\a_j(\var)$, 
$S=S_j$ be as in Remark 29. Let $\G_\var=\G_{j,\var}$ be the correspondent separatrix family from the previous Remark. Let $q$ be the vector function from Theorem 8 correspondent to the nonperturbed field and the sector $S$. There exist an $r>0$ and a family 
$\ov$ of domains in the $t$- plane containing $\a(\var)$ with the following properties:

1) the connected component of the intersection $\ov\cap (S^r\setminus\bigcup_{\a_s\neq\a}[0,\a_s(\var)])$ containing $\a(\var)$ tends to $S^r$, as $\var\to0$ (see footnote 1). 

2) For all $\var\neq0$ small enough the separatrix $\G_\var$ contains the graph $y=q_\var(t)$ of a vector function $q_\var(t)$ holomorphic  in $t\in\ov$ depending on the parameter, $q_\var(\a(\var))=0$,  such that 
$$\lim_{\var\to0}q_\var|_{\ov\cap S^r}=q\ \ \ \text{(see footnote 2).}$$
\endproclaim

In the case, when $n=2$, Theorem 10 is equivalent to Lemma 2 from Subsection 3.1. In higher dimensions it is proved in Subsections 4.2-4.3.

\subhead 4.2. Scheme of the proofs of Theorems 8 and 10 \endsubhead

For the proof of Theorems 8 and 10 we show (in the next Subsection) that in the conditions of Theorem 10 the separatrices $\G_\var$ contain graphs of vector-functions $q_\var(t)$ holomorphic in domains $\ov$ satisfying statement 1) of Theorem 10 such that inequalities (3) from Subsection 3.2 hold. (In particular, the family $\qe$ is normal in the sector $S^r$.) Then each limit $q$ of convergent sequence $q_{\var_m}$, $\var_m\to0$, satisfies the statements of the first part of Theorem 8: it is continuous in $\overline{S^r}$, vanishes at 0 and has bounded derivative; its graph over $S^r$ is tangent to the nonperturbed field. For the proof of Theorem 8 we consider a nondegenerate deformation of the given vector field (11) such that some its continuous singularity family corresponds to the sector $S$.    We show (in Subsection 4.4) that any vector function $q$ satisfying the three previous statements is $C^{\infty}(0)$ and has the asymptotic Taylor series $\hat q$ at 0. We prove uniqueness of such a vector function at the end of the paper. This will prove Theorem 8. 

The convergence statement of Theorem 10 follows from the statements on the vector-functions $\qe$ from the beginning of the previous item (inequalities (3) and statement 1) of Theorem 10) and the uniqueness statement of Theorem 8, as in the proof of Lemma 2.

\subhead 4.3. Proof of Theorem 10\endsubhead

 In higher dimensions a version of Theorem 10 was implicitly proved in \cite{6}  under the  assumption that $(11)_\var$ is tangent to the spaces $t=\a_i(\var)$ for all $\var$ and $i=0,\dots,k$ (i.e., in the notations of the item following Theorem 9 without loss of generality one can consider that $F\equiv p$ ($O_2\equiv0$)). 
Though a proof in \cite{6} was presented for vector field families obtained by projectivization of linear differential equations, it remains valid in this more general case. In difference with the two-dimensional case, in  general, in higher dimensions a vector field (11) does not necessarily have 
an integral hypersurface tangent to the $y$- space at 0, and in this case no its  deformation $(11)_\var$ is  orbitally analytically equivalent to a one with $F\equiv p$. Below we modify the proof from \cite{6} in order to 
make it valid for the general case. 

Without loss of generality we consider that the $t$- axis is tangent to the separatrix $G_\var$ of the field $(11)_\var$ at the singular point $(0,\a(\var))$. One can achieve this by applying appropriate continuous family of changes of the coordinates $y$. Then the separatrix $\G_\var$ contains the graph $y=\qe(t)$ of a vector function $\qe$ holomorphic in a neighborhood of $\a(\var)$ depending on $\var$, $\qe(\a(\var))=\qe'(\a(\var))=0$. For the proof of Theorem 10 it suffices to show that the vector  functions $q_\var$ are holomorphic in appropriate domains $\ov$ (satisfying statement 1) of Theorem 10) and satisfy inequalities (3) from Subsection 3.2. Let $K=(|\dot t|>|\dot y|)$ and $\widetilde K=\{|t-\a(\var)|>|y|\}$ be respectively the tangent cone field and  the correspondent cone  from the same Subsection. Inequalities (3) are equivalent to inclusions (4) from the same Subsection in the graphs $y=\qe(t)$. They hold a priori in a neighborhood of $\a(\var)$ depending on $\var$. To show that they hold in a large domain, we consider a pair of vector field families $(11)_{1,\var}=e^{i\theta_1}(11)_\var$, $(11)_{2,\var}=e^{i\theta_2}(11)_\var$,  $\theta_j\in\Bbb R$, which are constant multiples of $(11)_\var$, that possess the following properties: 

\item {1)}  The eigenvalue families of the linearization operators of the 
vector fields $(11)_{1,\var}$, $(11)_{2,\var}$ at the singular point $(0,\alpha(\var))$ correspondent to the eigenline tangent to the $t$- axis have positive
real parts, and their arguments  are bounded away from $\pi\Bbb Z$ uniformly in $\var$ small enough.
\smallskip
\item{2)}  The variables $y_i$ are split into two groups:  the
``stable'' $y_s=(y_{s_1},\ldots,y_{s_l})$ and the ``unstable''
$y_u=(y_{u_1},\ldots,y_{u_q})$ (one of them may be empty).  The tangent spaces of the $y_s$- and $y_u$- subspaces at $(0,0)$ are invariant with respect to the linearization operator of the nonperturbed field $(11)_0$ (and hence, to that of each of the nonperturbed fields $(11)_{1,0}$ and $(11)_{2,0}$ from the new families). The eigenvalues of the linearization operator of the field $(11)_{1,0}$ correspondent to the $y_s$- subspace have negative real parts; those correspondent to the 
$y_u$- subspace have positive real parts. Those of the linearization operator of the field $(11)_{2,0}$ correspondent to the $y_s$- ($y_u$-) subspace have oppositely positive (respectively, negative) real parts.

\

For a pair of real numbers $\theta_1$ and $\theta_2$ let us define the vector field families $w_j(\var):\ \dot t=e^{i\theta_j}p(t,\var)$ in the $t$- plane, $j=1,2$, which are the projection images to the latter of the restrictions to it of the vector field families $(11)_{j,\var}$.   The possibility of choice of $\theta_j$ satisfying the conditions 1) and 2) is implied by the following 

\remark{Remark 31} In the conditions of Theorem 10 let $SI\subset S$ be the 
closure in $S$ of the maximal subsector bounded by imaginary dividing rays. In the case, when $k=1$, let $V_\var=V_{j,\var}$ be the family of good sectors from Remark 29 correspondent to $\a(\var)$, $V$ be a good sector contained in the 
interior of the intersection $\cap_{\var}V_{\var}$.  In the case, when $k\geq2$, let $V=V_j$ be the correspondent sector from Remark 29. Let $\theta_1$, $\theta_2$ be chosen so that the correspondent nonperturbed vector fields 
$w_1(0)$, $w_2(0)$ of the families from the previous item have repelling rays lying in $V$ on different sides from $SI$: that correspondent to the former lies on the right 
and that correspondent to the latter lies on the left. Then $\theta_1$ and $\theta_2$ satisfy the previous conditions 1) and 2) for appropriate coordinates $(y_s,y_u)$. The proof of this statement is implicitly contained in remark 24 of \cite{6}. One can choose $\theta_i$ as at the beginning of the item so that $\theta_1>\theta_2$ and $\theta_1-\theta_2<\pi$. This follows from the fact that the angle of the sector $V$ is not greater than $\frac{\pi}k$ by definition.  
\endremark
\medskip
\noindent
\remark{ Remark 32}  In the conditions of Theorem 10 let $\theta_1,\theta_2,$ $(11)_{1,\var}$, 
$(11)_{2,\var}$, $(y_s,y_u)$ satisfy conditions 1) and 2) preceding Remark 31, 
$w_l(\var)$, $l=1,2$, be the correspondent vector field families in the $t$- line from the item preceding the previous Remark. Then $(0,\alpha(\var))$ is a hyperbolic
singular point for both $(11)_{1,\var}$ and $(11)_{2,\var}$.  The correspondent unstable manifold 
of the field $(11)_{1,\var}$ has tangent space at $(0,\a(\var))$ converging to  that  of the coordinate $(y_u,t)$ -space, as $\var\to0$. The tangent space at $(0,\a(\var))$ of the unstable manifold of the field $(11)_{2,\var}$ tends to that of the $(y_s,t)$-space.  In
particular, for small $\var$ the former intersects the latter transversally at $(0,\alpha(\var))$. The intersection of the unstable manifolds contains $\G_\var$ and coincides with it locally in a neighborhood of $(0,\a(\var))$. The conjoint singular point $\a(\var)$ of both perturbed fields $w_l(\var)$, $l=1,2$, is repelling, and the arguments of the correspondent multipliers are bounded away from $\pi\Bbb Z$ uniformly in $\var$. 
\endremark
\medskip
Let $\theta_1$ and $\theta_2$ be as in Remark 31. Recall that by assumption, 
the linearization matrices of the fields $(11)_{1,0}$ and $(11)_{2,0}$ at 0 are 
upper-triangular and block-diagonal (the eigenvalues of each block coincide; in particular, each block corresponds to either the $t$- axis, or a subspace of 
the coordinate $y_s$- ($y_u$-) space).  Without loss of generality we consider that their nondiagonal terms are small enough with respect to the nonzero eigenvalue real parts so that the sum of the moduli of the formers multiplied by 3 is less than the minimum of the moduli of the real parts of the nonzero eigenvalues. For the proof of inclusions (4) we consider the two following tangent cone fields:
$$
K_1=(|\dot t|+|\dot y_u|>3\,|\dot
y_s|)\,\,\text{and}\,\,K_2=(|\dot t|+|\dot y_s|>3|\dot y_u|),
$$
where
$$
|\dot y_u|=\sum\limits_i\,|\dot y_{u_i}|,\,\,\,|\dot
y_s|=\sum\limits_j\,|\dot y_{s_j}|.
$$
We also consider the correspondent cones 
$$
\widetilde K_1=\{|t-\alpha(\var)|+|y_u|>3\,|y_s|\},\,\,\widetilde
K_2=\{|t-\alpha(\var)|+|y_s|>3\,|y_u|\}.
$$
We show that
$$
T\,\Gamma_\var\subset K_1,\,\,T\,\Gamma_\var\subset K_2,\,\,\Gamma_\var
\subset\widetilde K_1\cap\widetilde K_2\tag12
$$
over appropriate set $\Omega_\var$. Then (12) will imply (4), since 
$K\supset(K_1\cap K_2)$ and $\widetilde K\supset\widetilde K_1\cap\widetilde K_2$.
\medskip
 The fields 
$(11)_{1,\var}$ and $(11)_{2,\var}$ are tangent to $\G_\var$. Their restrictions
to the latter have repelling singularity at $(0,\a(\var))$. 
For all $\var$ small enough inclusions (12) hold locally in a neighborhood of the singularity (that depends on $\var$). To show that (12) hold in a
large domain, we use the fact that the cone field $K_1$ ($K_2$) is invariant with
respect to the vector field $(11)_{1,\var}$ ($(11)_{2,\var}$) in a
neighborhood $U$ of zero in the phase space independent on $\var$  for
all $\var$ small enough. The proof of this statement uses only conditions 1), 2) preceding Remark 31 and the assumption from the item preceding the definition of the cone field $K_1$ and is analogous to that of proposition 3 in  \cite{6}.  
 
Without loss of generality we consider that  $U=U_y\times U_t$, $U_y=\{|y_s|<2\delta\}\times\{|y_u|<2\delta\}$, $U_t=
\{|t|<\delta\}$. Let $L\ov^1$ and  $L\ov^2$  be the subdomains of $\G_\var\cap U$ saturated by the  trajectories of the field $(11)_{1,\var}$ (respectively, $(11)_{2,\var}$) in 
$\G_\var\cap U$ that go from $(0,\a(\var))$. Let $L\ov$ be the connected component of the intersection $L\ov^1\cap L\ov^2$ that contains $(0,\a(\var))$. Inclusions (12) hold in $L\ov$ by invariance of the cone fields and the fact that they hold locally in a neighborhood of the  singularity. For the proof of Theorem 10 it suffices to show that (12) hold in the graphs $y=q_\var(t)$ of vector functions $q_\var$ holomorphic in appropriate domains $\ov$ satisfying statement 1) of Theorem 10. To do this, it suffices to show that $L\ov$ contain such graphs. We prove this for domain families $\ov$ of the following type.  

\definition{Definition of the domains $\ov$} Let  $U_t=\{|t|<\delta\}$ be a neighborhood of zero in the $t$- plane, $\theta_i$, $w_i(\var)$, $i=1,2$, be as in Remark 31. Let $\ov^i$ be the family of domains saturated by the trajectories of the field $w_i(\var)$ in $U_t$ that go from $\a(\var)$. Define $\ov$ to be the connected component of the intersection $\ov^1\cap\ov^2$ that contains $\a(\var)$ (see Fig.4a in the case, when $k=1$).
\enddefinition
\define\wto{L\ov'}

There exist $\theta_1$ and $\theta_2$ as in Remark 31 such that the correspondent family $\ov$ from the previous Definition satisfies statement 1) of Theorem 10. Namely, this is the case, provided that the correspondent repelling rays from Remark 31 are close enough to the imaginary dividing rays 
that form the boundary of the sector $SI$ from the latter, as in the proof of 
the analogous statement in Subsection 3.2.

For the proof of Theorem 10 it would be sufficient to show that $L\ov$ contains the graph of a vector function $q_\var(t)$ holomorphic in $\ov$. In the case, when $O_2\equiv0$, $L\ov$ is such a graph itself: it is  1-to-1 projected onto $\ov$ (see subsection 5.B in \cite{6}).  This is the place in the proof of the particular case of Theorem 10 from \cite{6} we used the last equality. 

Let $\ov$ be a given domain family as in the previous Definition. In the case,  when $O_2\not\equiv0$, for the proof of Theorem 10 we 
consider another pair of numbers $\theta_i'$, $i=1,2$, satisfying the conditions of Remark 31 such that $\theta_1>\theta_1'>\theta_2'>\theta_2$: the ordering of the "arguments" $\theta_2, \theta_2',\theta_1',\theta_1$ defines  the counterclockwise order of the correspondent radial rays. We consider the correspondent domain family $\wto$ from the item preceding the previous Definition.  We show that the new domain 
$\wto$ (where inclusions (12) hold by the discussion from the same place) contains the graph of a vector function $q_\var(t)$ holomorphic in $\ov$, whenever $\delta$ and $\var$ are small enough. This will prove Theorem 10. 

The projection of the new domain 
$\wto$ to the $t$- plane is a local diffeomorphism. This follows from the fact that the tangent lines to $\wto$ are contained in the tangent cone field $K\supset K_1\cap K_2$. Thus, the inverse map $\ov\to\wto$ is well-defined a priori in a neighborhood of the point $\a(\var)$ and extends analytically to  
any its simply connected neighborhood disjoint from the projection image of 
$\partial\wto$. We show that the whole domain $\ov$ is disjoint from the latter, whenever $\delta$ and $\var$ are small enough. This will imply the last statement from the last item.  To do this, we use the following properties of $\wto$. 

\remark{Remark 33} Let $\theta_1$, $\theta_2$ be as in Remark 31, $U$, $U_y$,  $U_t$,  $L\ov$ be as in the two items preceding the previous Definition. The domain $L\ov$ meets the boundary of the neighborhood $U$ at points of the cylinder $U_y\times\partial U_t$. The domain $L\ov$ is invariant with respect to the vector fields $-(11)_{1,\var}$ and $-(11)_{2,\var}$, i.e., for any $a\in L\ov$ the arcs of 
their trajectories that connect $a$ to the singular point $(0,\a(\var))$ are 
contained in $L\ov$. The first statement of the item follows from the inclusion 
$L\ov\subset \widetilde K$ ($TL\ov\subset K$) in 
the same way, as in the proof of lemma 1 in subsection 5.B of \cite{6}. Let us prove the second statement of the item, e.g., for the first vector field. Let $L\ov^1$, $L\ov^2$ be as in the item preceding the previous Definition. The domain 
$L\ov^1$ ($L\ov^2$) is $-(11)_{1,\var}$- (respectively, $-(11)_{2,\var}$-) invariant. Therefore, by the latter and definition, the boundary of the domain $L\ov$ consists of parts of $\partial L\ov^1$ and 
some (may be, semiinfinite) arcs of trajectories of the field $(11)_{2,\var}$ that start and end at 
$\partial L\ov^1$. Each of these arcs splits the domain $L\ov^1$ into two 
connected components. Let $L$ be one of these arcs. By $L\ov^0$ denote the 
 correspondent splitting component of $L\ov^1$ that contains the singular point. It suffices to show that the field $(11)_{1,\var}$ 
is directed outside $L\ov^0$ at $L$. This follows from the fact that 
the angle between the vector fields is constant (in particular, $(11)_{1,\var}$ has constant angle with $L$) and there exists a point of $L$ where $(11)_{1,\var}$ is directed outside $L\ov^0$. The last statement is implied 
by the definition of the domain $L\ov^1$: each point of the complement $L\ov^1\setminus\L\ov^0$ is connected to the singular point $(0,\a(\var))$ by trajectory of the field $(11)_{1,\var}$ in $L\ov^1$. 
\endremark

\proclaim{Corollary 6} In the conditions of the previous Remark the boundary 
$\partial L\ov$ consists of the two following parts: 

1) arcs in the cylinder $U_y\times\partial U_t$ where both vector fields $(11)_{1,\var}$ and $(11)_{2,\var}$ are directed outside the latter;

2) arcs of semitrajectories of the fields $(11)_{1,\var}$ and $(11)_{2,\var}$ in $\overline U$ that start at some points of their tangency with the cylinder.
\endproclaim

\example{Example 7} Let $k=1$. Then in the conditions of the previous Remark   one can show that for $U$ and $\var$ small enough the part 1) 
of the boundary $\partial L\ov$ consists of a single arc and 
the part 2) consists generally of two positive semitrajectories of the fields $(11)_{1,\var}$ and $(11)_{2,\var}$ that converge to the other singular 
point $(0,-\a(\var))$. (In some exceptional cases it consists of two finite arcs of these semitrajectories that have a common end.) The projection images of these semitrajectories in the $t$- plane are close to the semitrajectories of the correspondent vector fields $w_1(\var)$ and $w_2(\var)$ that form the boundary of $\ov$ (see Fig. 4a and Proposition 6 below). The semitrajectories of the latters intersect each other, so, the correspondent part of the boundary of the domain $\ov$ consists of their arcs having a common end. In particular, the projection images of the previous  semitrajectories (forming $\partial L\ov$) intersect each other. At the same time, one can show that generically these semitrajectories do not intersect each other (so, the projection images in the $t$- plane of some disjoint parts of $\wto$ overlap). The proof of this statement is omitted to save the space. 
\endexample 

Now let us show that the projection image in the $t$- plane of the boundary 
$\partial\wto$ does not intersect $\ov$, whenever $U$ and $\var$ are 
small enough. The projection image of the part 1) of the boundary from Corollary 6 does not intersect  $\ov$ by definition. Thus, it suffices to show that no projection of arc from the part 2)  
(which starts outside $\ov$) can enter $\ov$.  To do this, we use the      following 

\proclaim{Proposition 5} Let $p(t,\var)=\pr$ be a monic polynomial vector field family that satisfies the asymptotic singularity polygon regularity statement from Remark 10 in the case, when $k\geq2$. Let ($\theta_1$, $\theta_2$), ($\theta_1',\theta_2')$, be real number pairs such that $\theta_1>\theta'_1>\theta_2'>\theta_2$, $\theta_1-\theta_2<\pi$ (i.e., the ordering of the "arguments" $\theta_2, \theta_2', \theta_1',\theta_1$ defines the counterclockwise order of the correspondent radial rays in complex plane). Let ($w_1(\var)$, $w_2(\var)$), 
($w_1'(\var)$, $w_2'(\var)$) be the correspondent vector field family pairs 
from the item preceding Remark 31. Let $\a(\var)=\a_i(\var)$ be their conjoint  singular point family such that the correspondent multipliers of the four perturbed fields have positive real parts and their arguments are bounded away from $\pi\Bbb Z$. Let $U_t$ be a neighborhood of zero in the $t$- line, $\ov$ be the domain family from the previous Definition correspondent to the first vector field family pair. The domain $\ov$ is invariant with respect to the fields 
$-w_1'(\var)$ and $-w_2'(\var)$. Each of the latters is directed strictly 
inside $\ov$, and moreover its angle with the arcs of the boundary $\partial\ov$ is bounded away from $\pi\Bbb Z$ uniformly in $U_t$ and $\var$. 
\endproclaim

\demo{Proof} The domain $\ov$ is $-w_1(\var)$- and 
$-w_2(\var)$- invariant, as is $L\ov$ with respect to the correspondent vector fields in Remark 33. Let $t\in\partial\ov$. Let us show that the fields $w_1'(\var)$ and $w_2'(\var)$ are directed outside $\ov$ at $t$. Let $W=W(t,\var)$ be the
sector in the tangent plane at $t$ with angle less than $\pi$ formed by the 
vectors of the fields $w_1(\var)$ and $w_2(\var)$. Each radial vector in $W$ is directed outside $\ov$. This follows from the statement that  both $w_1(\var)$ and $w_2(\var)$ are not directed inside $\ov$ in its boundary (the previous invariance statement). The vectors of the fields $w_1'(\var)$ and $w_2'(\var)$ are contained in $W$ and have constant angles with its boundary rays by definition.  This proves Proposition 5.\enddemo

\proclaim{Corollary 7} In the conditions of Proposition 5 the domain $\ov'$ from the previous Definition correspondent to the fields $w_i'(\var)$, $i=1,2$,  contains $\ov$. \endproclaim

Let $(11)_{1,\var}'$, $(11)_{2,\var}'$ be the vector fields from the beginning of the proof of Theorem 10 correspondent to the new numbers $\theta_1'$ and $\theta_2'$, $\wto$ be the correspondent domain from the item preceding the previous Definition. The arcs from the part 2) of $\partial\wto$ lie in trajectories of the fields $(11)_{1,\var}'$ and $(11)_{2,\var}'$. 
Now the last statement preceding Proposition 5 will follow from the latter and 
the fact that the projection images in the $t$- plane of the vectors of these fields restricted to the closure $\overline{L\ov'}$ have arbitrarily small angles with the  
field $w_1'(\var)$ (respectively, $w_2'(\var)$), provided that $U$ and $\var$ are small enough. Then these angles can be made so small that these projection vectors  could not be directed inside $\ov$, as the fields $w_i'(\var)$ in  Proposition 5. Thus, no projection image of arc from part 2) of the boundary $\partial L\ov'$, which starts outside $\ov$, can enter the latter.

 The projection to the $t$- line of the closure $\overline{L\ov'}$ is local diffeomorphism outside the singular points of the field $(11)_\var$: the restriction of $(11)_\var$ to $\overline{L\ov'}$ lies in the cone $K$ closure field, since so does the restriction to $L\ov'$. The angles from the last item are equal to the angle between the vector field $(11)_{\var}$ and the lifting to $\overline{\wto}$  of the field $p(t,\var)$. The angle bound statement from the last item is implied by the following 
  
\proclaim{Proposition 6} Let $(11)_\var$, $\a$ be as in Theorem 10. Let $K$  and $\widetilde K$ be the cone field and the cone from the beginning of the Subsection. 
 For any $\var$ (including $\var=0$) let  $K_\var\subset \overline{\widetilde K}\setminus\{(0,\a_i(\var)); i=0,\dots,k\}$ be the subset of points  
where the vectors of the field $(11)_\var$ are contained in the cone closure field $\overline K$. 
  Let $p(t,\var)$ be the polynomial family from the item following Theorem 9, $(11)''_\var$ be the vector field in $K_\var$ contained in the complex tangent 
line field correspondent to $(11)_\var$ and projected to the vector 
field $\dot t=p(t,\var)$ in the $t$- plane. Then for any $\sigma>0$ there 
exists a neighborhood $U$ of zero in the phase space such that for all $\var$ 
small enough the inequality $|(11)_\var-(11)''_{\var}|<\sigma|(11)''_\var|$  holds in $K_\var\cap U$.
\endproclaim
 
Proposition 6 is proved below.

In the conditions of the item following Corollary 7 for any $\sigma>0$ 
and any $U$, $\var$ small enough (dependently on $\sigma$) the angle between 
the field $(11)_{l,\var}'|_{\wto}$ and the lifting to $\wto$ of the field $w_i'(\var)$, $i=1,2$, is less than $\sigma$ everywhere in $\overline{\wto}$. This follows from Proposition 6 and the inclusion $L\ov'\subset K_\var$, which holds whenever $U$ and $\var$ are small enough. This together with Proposition 5 and Corollary 6 proves that for $U$ and $\var$ small enough the domain $\ov$ is disjoint from the projection image of the boundary of $L\ov'$. Theorem 10 is proved modulo Proposition 6.

\demo{Proof of Proposition 6} Let $B$ be the matrix from expression (11) for the   nonperturbed vector field, $G(y,t,\var)$, $F(y,t,\var)$,  $p(t,\var)$  
be the (vector) functions from the definition of $(11)_\var$ and the item following Theorem 9. For the proof of Proposition 6 it suffices to show that $F(y,t,\var)=p(t,\var)(1+o(1))$, as $(y,t,\var)$ tends to 0 so that $(y,t)\in K_\var$. Suppose the contrary, i.e.,  there exists a $c>0$ such that the inequality 
$|F-p|>c|p|$ holds along a sequence $(y_m,t_m,\var_m)\to0$, $(y_m,t_m)\in K_{\var_m}$. 
 By the assumptions from the item following Theorem 9, $F(y,t,\var)=p(t,\var)+o(|y|)$. Therefore, both $p$ and $F$ are $o(|y|)$ along this sequence. Hence, so is $G$, which follows from the  same statement for $F$ and the inequality $|G|\leq|F|$ (the definition of the set $K_\var$ in Proposition 6). On the other hand, 
$G(y,t,\var)=By+o(|y|)$ along the same sequence, which follows immediately from the expression for $G$ in the item following Theorem 9 and 
the asymptotic $p(t,\var)=o(|y|)$ along the sequence proved before. Therefore, $G(y,t,\var)\neq o(|y|)$ along the sequence (the matrix $B$ is nondegenerate).   This contradicts 
the previous statement that $G(y,t,\var)=o(|y|)$ along the sequence. The obtained contradiction proves Proposition 6. The proof of Theorem 10 is completed. 
\enddemo

\subhead 4.4. Proof of Theorem 8 \endsubhead

Let (11), $S$ be as in Theorem 8. For the proof of 
 Theorem 8 we consider a nondegenerate deformation $(11)_\var$ of the field (11) with the following properties.

\remark{Remark 34} For any vector field (11) and sector $S$ as in Theorem 8 there exists a nondegenerate deformation $(11)_\var$ of (11) (satisfying the 
assumptions of the item following Theorem 9) with a singularity coordinate  family $\a(\var)$ such that the correspondent sector from Remark 29 contains the same  imaginary dividing rays, as $S$. Indeed, the right-hand side in (11) is $O(|y|)+t^{k+1}g(t)$ where $g$ is a holomorphic vector function. Changing the multiplier $t^{k+1}$ in the last formula to any given degree $k+1$ monic polynomial family $p(t,\var)$, $p(t,0)=t^{k+1}$, yields a deformation $(11)_\var$ satisfying the assumptions from the item following Theorem 9. Let $\theta\in\Bbb R$ be a fixed constant. Let us choose a polynomial family $p(t,\var)$ dependent on the real parameter $\var$ as follows:  $p(t,\var)=t^{k+1}-e^{i\theta}\var=\prod_{j=0}^k
(t-e^{i(\theta+\frac{2\pi}{k+1}j)}\var^{\frac1{k+1}})$, so, the correspondent root families $\a_j(\var)=e^{i(\theta+\frac{2\pi}{k+1}j)}(\var)^{\frac1{k+1}}$  form a regular polygon family centered at 0 with constant vertex radial rays. Let $SI\subset S$ be the maximal subsector bounded by imaginary dividing rays. A deformation $(11)_\var$ defined by a polynomial family $p(t,\var)$ as above will be a one we are looking for, provided that $\theta$ is chosen so that the radial ray 
of some root family $\a=\a_j(\var)$ is close enough to the bissectrix of the sector $SI$ and no root radial ray coincides with a real dividing ray.  
\endremark

 Let $\G_\var$ be the separatrices of the perturbed fields $(11)_\var$ at $(0,\a(\var))$ from Remark 30. The separatrices $\G_\var$ contain the graphs $y=q_\var(t)$ of holomorphic vector functions $q_\var(t)$ defined in domains $\ov$ satisfying statement 1) of Theorem 10 such that $|q_\var'|<1$, $|q_\var(t)|<|t-\a(\var)|$ (for $t\neq\a(\var)$). These statements were proved in the previous Subsection. Then the vector function family $q_\var$ is normal, so any sequence $q_{\var_m}$, $\var_m\to0$, contains a subsequence convergent uniformly in compact subsets in $S^r$. Each limit $q(t)$ of convergent  sequence $q_{\var_m}$, $\var_m\to0$, as $m\to\infty$, is a vector function $q(t)$ holomorphic in the sector $S^r$ and continuous in its closure such that a) $q(0)=0$; b) $|q'|_{S^r}\leq1$; c) the graph $y=q(t)$ over the sector $S^r$ is contained in a phase curve of the nonperturbed field. (Statements a) and b) follow from the previous inequalities.) Let us show that a vector function $q(t)$ with properties a), b), c) satisfies the statements of Theorem 8.

Firstly let us prove that $q$ has asymptotic Taylor series at 0 (coinciding with the formal central manifold series $\hat q$ from the beginning of the Section). To do this, let us prove the correspondent asymptotic Taylor formula. Namely, for $m\in\Bbb N$ define $q_m(t)$ to be the Taylor polynomial of degree at most $m$ coinciding with the $m$- th partial sum of the series $\hat q$. Let us show that $q(t)-q_m(t)=o(t^{m})$, as $t\to0$ for any $m$.  
Without loss of generality we consider that $q_m\equiv0$. One can achieve this by applying the change $y\mapsto y-q_m(t)$ of the coordinates $y$. Then the vector field (11) takes the form 
$$\cases
&\dot y=G(y,t)=By+O(|y|^2+|y||t|)+ O(t^{m+1})\\
&\dot t=F(y,t)=t^{k+1}+ O(|y|^2).\endcases\tag13$$
Let us show that $q(t)=o(t^{m})$, as $t\to0$. In the proof of this statement we use the asymptotic formula  
$$F(y,t)|_{y=q(t)}=t^{k+1}(1+o(1)), \text{as}\  t\to0.\tag14$$ 
(This is equivalent to the statement that the angle between the field (11) and the lifting to the graph $y=q(t)$ of the field $\dot t=t^{k+1}$ in the $t$- line tends to 0, as $t\to0$, cf. the proof of Proposition 6). This follows from inequality b) in the last item (which implies that the graph $y=q(t)$ is contained in the closure of the cone $\widetilde K$ from the beginning of the previous Subsection) and Proposition 6  applied to the nonperturbed field and its singular point 0. Consider the coordinates $(y_s,y_u)$ and the vector fields  $(11)_{1,0}$,  $(11)_{2,0}$, $w_1(0)$, $w_2(0)$ from Remark 31 chosen so that the two latters have   repelling rays in $S$. (In general, the sector $S$ from Theorem 8 does not contain the whole sector $V$ from Remark 31. A priori, it contains the closure of the sector $SI$ from the latter. One can choose the two last fields as in the same Remark so that the correspondent repelling rays will be close enough to the boundary of the sector $SI$ so that they fit  the sector $S$.) The angle between $(11)_{i,0}$ and the lifting to the graph $y=q(t)$ of the field $w_i(0)$ tends to 0, as $t\to0$. Then the trajectories of the field $-(11)_{i,0}$ that start in the graph over the sector $S^{r_0}$ with appropriate $r_0>0$ converge to 0 so that their projections to the $t$- plane do not leave the sector $S^r$ and converge to zero with the asymptotic tangency to the repelling  ray of the field $w_i(0)$. Recall that the linearization operators of the fields $-(11)_{1,0}$, $-(11)_{2,0}$ have invariant subspaces tangent to the $y_s$- and $y_u$- spaces; the eigenvalues of that of the field $-(11)_{1,0}$ in the first subspace have positive real parts, and those in the second subspace have 
negative real parts; those correspondent to the second field have oppositely 
negative (respectively, positive) real parts. The nondiagonal terms of the 
linearization matrices of these fields are small with respect to the real parts 
of nonzero eigenvalues, more precisely, satisfy the assumptions from the 
item following Remark 32. 

We prove the asymptotic formula $q(t)=o(t^m)$ by contradiction. Suppose the contrary. This means that there exists a $c>0$ such that there exists arbitrarily small $t_0\in S^r$ where 
$|q(t_0)|>c|t_0|^m$. We consider that $|q_s(t_0)|\geq|q_u(t_0)|$, so $|q_s(t_0)|>\frac c2|t_0|^m$ (the opposite case is analyzed analogously). Let us show that the norm $|q_s(t)|$ increases along the trajectory of the field $-(11)_{1,0}$ that starts at $(q(t_0), t_0)$, provided that $t_0$ is close enough to 0. Then 
this will contradict the convergence of this trajectory to 0 and  will prove the asymptotic Taylor formula. For the proof of the increasing of this norm we consider the subdomain $K_m=\{ |y_s|>|y_u|,\ |y_s|>\frac c2|t|^m\}\cap\{y=q(t)\}$ of the graph of the vector function $q(t)$. The closure of the domain $K_m$ contains the point $(q(t_0), t_0)$, by definition. We use the fact that $K_m$ is $-(11)_{1,0}$- invariant over a neighborhood of 0 small enough in $S^r$ (more precisely, for any point of $\partial K_m$ close enough to 0 the correspondent positive semitrajectory of the field $-(11)_{1,0}$ 
enters $K_m$ locally). This statement follows from the definition of the coordinates $(y_s,y_u)$, (13),  (14) and is proved by a straightforward calculation of the derivative of the ratios $\frac{|y_s|}{|y_u|}$, $\frac{|y_s|}{|t|^m}$, along the field 
$-(11)_{1,0}$ in $\partial K_m$: these derivatives will be positive, whenever $(y=q(t),t)\in\partial K_m$ is close enough to 0, as in the proof of proposition 3 from \cite{6}. For the proof of the previous $|q_s|$ increasing statement it suffices to show that the norm $|y_s|$ increases along the trajectories of the field $-(11)_{1,0}$ in $\overline{K_m}$ at all the points of the latter close enough to 0. Indeed, the increment of the norm $|y_s|$ along a small trajectory segment consists of the following parts: the contribution of the linear terms $By$ in the right-hand side of  (13) and that of the nonlinear terms. The former is positive, and moreover it is greater than $|y_s|$ times some positive constant independent on $(y,t)\in K_m$ (by definition of the coordinate subspaces $y_s$ and $y_u$). It dominates the nonlinear term contribution, whenever $(y,t)\in \overline{K_m}$ is close enough to 0. This follows from  the fact that the nonlinear terms are $o(|y|)+o(|t|^{m+1})=o(|y_s|)$, as $(y,t)\to0$ so that $(y,t)\in K_m$ (the definition of $K_m$). This proves the asymptotic Taylor formula for $q(t)$. The statement of  Theorem 8 that $\hat q$ is the asymptotic Taylor series of the function $q$ is proved modulo continuity of the derivatives of $q$ at 0.    
In the proof of the last statement we use the asymptotic Taylor 
formula proved before and the differential equation for the vector function $q$:  $$\frac{dq}{dt}=\frac{G(y,t)}{F(y,t)}|_{y=q(t)}.$$ 
Let us prove continuity of $q'$. For higher derivatives the proof is analogous.  Let us show that $q'(t)\to0$, as $t\to0$. 
Let $q_{k+1}(t)$ be the Taylor polynomial from the previous asymptotic Taylor formula. Without loss of generality we consider that $q_{k+1}\equiv0$, as in the item preceding (13). Then $q(t)=o(t^{k+1})$ by the asymptotic Taylor formula, hence $G(q(t),t)=Bq(t)+O(|q(t)|^2+|tq(t)|+|t|^{k+2})=o(t^{k+1})$. This together with the previous differential equation and (14) implies that $q'(t)=o(1)$, as $t\to0$. This proves the continuity of the extension of $q'$ to 0 by zero. 

Now let us prove the uniqueness of a vector function $q$ with properties a) and c) from the beginning of the proof of Theorem 8 and bounded derivative. Suppose there exists another such vector function $\tilde q\neq q$. Without loss of generality we consider that the module of its derivative is not greater than 1 
(i.e., inequality b) from the beginning of the Subsection holds). One can achieve this by applying appropriate linear change $y\mapsto\lambda y$ of the coordinates $y$. Then by the previous discussion, $\tilde q$ has asymptotic Taylor series at 0. We consider that $q\equiv0$. One can achieve this by applying the change $y\mapsto y-q(t)$ of the variable $y$ over the sector $S^r$. (Then the new vector field (11) will be holomorphic over $S^r$ and in general, only $C^{\infty}$ at 0 with asymptotic Taylor series.) Let us show that $\tilde q\equiv0$. 

Suppose the contrary: there exists a $t_0\in S^r$ where $\tilde q(t_0)\neq0$. We consider that $|\tilde q_s(t_0)|\geq|\tilde q_u(t_0)|$ (the opposite case is analyzed analogously). The semitrajectory of the field $-(11)_{1,0}$ that starts at 
$(\tilde q(t_0),t_0)$ converges to 0 so that its projection to the $t$- line does not leave the sector $S^r$, provided that $t_0$ is small enough, as in the proof of the asymptotic Taylor formula. Let us show that the norm $|\tilde q_s(t_0)|$ increases 
along this trajectory, if $t_0$ is small enough. This will contradict its  convergence to 0. For the proof of the increasing of this norm we consider the tangent cone field $(|\dot y_s|>|\dot y_u|)$ at the graph of the function $\tilde q$. We use the fact that this cone field is $-(11)_{1,0}$- invariant in appropriate neighborhood of zero in this graph (and hence, so is the correspondent domain $K_{s,u}=\{|y_s|>|y_u|\}\cap\{y=\tilde q(t)\}$, whose closure contains the point $(\tilde q(t_0),t_0)$, by definition). The proof of these statements is analogous to that of proposition 3 in \cite{6}. Now the proof of the $|\tilde q_s|$- increasing statement repeats that of the analogous statement from the proof of the asymptotic Taylor formula with the change of 
$K_m$ to $K_{s,u}$. This contradicts the convergence to 0 of the trajectory from the beginning of the item and proves that $\tilde q\equiv0$. Theorem 8 is proved.

\head Acknowledgements\endhead

I am grateful to Yu.S.Ilyashenko for the statement of the problem. 
I am grateful to him, to A.A.Bolibroukh and J.-P.Ramis for helpful discussions.
The paper was written when I was visiting Institute for Advanced Study 
(Princeton, USA) and Max-Planck Institut f\"ur Mathematik (Bonn, Germany). I wish to thank both Institutes for hospitality and support. 

\head References\endhead

1. Yu.S.Ilyashenko, Nonlinear Stokes Phenomena. - Adv.Soviet Math., 14, AMS, 
1993. 

2. J.Martinet, J.-P. Ramis. Probl\`emes de modules pour des \'equations 
diff\'erentielles non lin\'eaires du premier ordre. - Inst. Hautes \'Etudes Sci.  Publ. Math. No. 55 (1982), 63--164. 

3.  R.Garnier.\ \  
 Sur\ \  les singularit\'es \ irreguli\'eres \ \ des \'equations\ \ 
diff\'erentielles lin\'eaires.
\rm J. Math. Pures et Appl., 8$^e$ s\'erie, 2 (1919), pp. 99--198.

4. J.Martinet. Remarques sur la bifurcation n\oe ud-col dans le domaine 
complexe. - Asterisque No.150-151 (1987), 131-149.

5. C.Zhang. Quelques \'etudes en th\'eorie des \'equations fonctionnelles
et en analyse combinatoire. 
\rm Th\'ese 1994.  Institut de Recherche Mathematique Avanc\'ee,
Univ. Louis Pasteur et CNRS (URA 01).

6. A.A.Glutsuk. Stokes operators via limit monodromy of generic perturbation. - 
To appear in J. Dynamical and Control Systems. 

7. V.I.Arnold. Geometric methods in ordinary differential equations. - 
Springer-Verlag, 1983.  Translated 
from the Russian original "Dopolnitelnyie glavy obyknoviennyh differentsialnyh 
uravneniy"( "Nauka", Moscow, 1978).  

8. J.Ecalle. Invariants holomorphes simples des transformations de 
multiplicateur 1. - C. R. Acad. Sci. Paris S\'er. A-B 276 (1973), A375-A378.

9. S.M.Voronin. Analytic classification of germs of conformal mappings 
$(\Bbb C,0)\rightarrow(\Bbb C,0)$ (in Russian). - Funktsionalnyi 
Analiz i Prilozheniya, 15 (1981), No 1, 1-17.

10. L.Carleson, Th.Gamelin. \ Complex dynamics. - Universitext: Tracts in 
Mathematics. Springer-Verlag, New York, 1993.

11. J.-P. Ramis, Y.Sibuya. Hukuhara domains and fundamental existence and uniqueness theorem for asymptotic solutions of Gevrey type. - Asymptotic Anal. 2 (1989), no.1, 39-94.

12. V.I.Arnold, Yu.S.Ilyashenko. 
\sl ``Ordinary differential equations''.
\rm In the volume ``Dynamical Systems--1'' of the series ``Itogi nauki i
tekhniki.  Sovremennyie problemy matematiki.  Fundamentalnyie
napravlenia'', VINITI publisher, Moscow (1985).

\enddocument